\documentclass{amsart}

\usepackage{amsmath,amsthm,amsfonts,amscd,amssymb} 
\usepackage{verbatim}      	
\usepackage[normalem]{ulem}
\usepackage{epsfig}         	
\usepackage{url}		
\usepackage{epic,eepic,latexsym}
\usepackage{graphicx}
\usepackage{color}
\usepackage{lcd}
\usepackage{mathrsfs}
\usepackage[centering,text={15.5cm,22cm},
        marginparwidth=20mm,dvips]{geometry}
\usepackage[linktocpage]{hyperref}
\usepackage{lscape}
\usepackage{tabularx}
\usepackage{marginnote}
\usepackage[all]{xy}
\usepackage{enumitem}
\usepackage{stmaryrd}
\usepackage{upgreek}

\DeclareMathAlphabet{\mathpzc}{OT1}{pzc}{m}{it}

\DeclareMathOperator{\Hom}{Hom}

\DeclareMathOperator{\codim}{codim}

\DeclareMathOperator{\rank}{rank}

\DeclareMathOperator{\Spec}{Spec}
\DeclareMathOperator{\id}{Id}

\DeclareMathOperator{\NE}{NE}

\DeclareMathOperator{\Tr}{Tr}

\DeclareMathOperator{\val}{val}

\DeclareMathOperator{\TV}{\s{TV}}
\DeclareMathOperator{\sign}{sgn}

\DeclareMathOperator{\trop}{trop}
\DeclareMathOperator{\prin}{prin}
\DeclareMathOperator{\out}{out}
\DeclareMathOperator{\In}{in}
\DeclareMathOperator{\Mult}{Mult}
\DeclareMathOperator{\mult}{\wt{Mult}}

\DeclareMathOperator{\vir}{vir}
\DeclareMathOperator{\ev}{ev}
\DeclareMathOperator{\Aut}{Aut}
\DeclareMathOperator{\Ends}{Ends}

\DeclareMathOperator{\scat}{Scat}

\DeclareMathOperator{\Joints}{Joints}

\DeclareMathOperator{\Supp}{Supp}

\DeclareMathOperator{\ad}{ad}

\DeclareMathOperator{\vdim}{vdim}
\DeclareMathOperator{\N}{N}

\DeclareMathOperator{\GW}{GW}
\DeclareMathOperator{\Flags}{Flags}

\DeclareMathOperator{\pt}{pt}
\DeclareMathOperator{\uf}{uf}
\DeclareMathOperator{\Bl}{Bl}
\DeclareMathOperator{\SH}{SH}
\DeclareMathOperator{\QH}{QH}

\DeclareMathOperator{\topp}{top}
\DeclareMathOperator{\forget}{Forget}
\DeclareMathOperator{\naive}{naive}
\DeclareMathOperator{\KN}{KN}

\DeclareMathOperator{\irr}{irr}

\let\?=\overline
\let\:=\colon
\let\llb=\llbracket
\let\rrb=\rrbracket
\let\bb=\mathbb

\let\rar=\rightarrow
\let\f=\mathfrak
\let\s=\mathcal

\let\wh=\widehat
\let\wt=\widetilde

\def\risom{\buildrel\sim\over{\smashedlongrightarrow}}
 \def\smashedlongrightarrow{\setbox0=\hbox{$\longrightarrow$}\ht0=1.25pt\box0}

\newcommand {\AAA} {\vec{{\bf A}}}
\newcommand {\A} {{\bf A}}
\newcommand {\kk} {\Bbbk}

\newcommand {\ww} {{\bf w}}

\newcommand {\pp} {{\bf p}}

\newcommand {\SSS} {{\bf S}}

\theoremstyle{plain}
 \newtheorem{thm}{Theorem}[section]
 
 \newtheorem{lem}[thm]{Lemma}
  
  \newtheorem{prop}[thm]{Proposition}
  \newtheorem{conj}[thm]{Conjecture}
   \newtheorem{cor}[thm]{Corollary}

\theoremstyle{definition}
 \newtheorem{dfn}[thm]{Definition}
  
 \newtheorem{ntn}[thm]{Notation}

\theoremstyle{remark} 
 \newtheorem{rmk}[thm]{Remark}

    \newtheorem{asss}[thm]{Assumptions}

\newenvironment{myproof}[1][\proofname]{\proof[#1]\mbox{}}{\endproof}

\title{Theta bases and log Gromov-Witten invariants of cluster varieties}
\author{Travis Mandel}
\address{Department of Mathematics\\
University of Oklahoma\\
Norman, OK 73019\\
USA}
\email{tmandel{\char'100}ou.edu}
\thanks{The author was supported by the National Science Foundation RTG Grant DMS-1246989, and later by the Starter Grant ``Categorified Donaldson-Thomas Theory'' no. 759967 of the European Research Council.}

\begin{document}

\begin{abstract}
Using heuristics from mirror symmetry, combinations of Gross, Hacking, Keel, Kontsevich, and Siebert have given combinatorial constructions of canonical bases of ``theta functions'' on the coordinate rings of various log Calabi-Yau spaces, including cluster varieties.  We prove that the theta bases for cluster varieties are determined by certain descendant log Gromov-Witten invariants of the symplectic leaves of the mirror/Langlands dual cluster variety, as predicted in the Frobenius structure conjecture of Gross-Hacking-Keel.  We further show that these Gromov-Witten counts are often given by naive counts of rational curves satisfying certain geometric conditions.  As a key new technical tool, we introduce the notion of ``contractible'' tropical curves when showing that the relevant log curves are torically transverse. \end{abstract}

\maketitle

\setcounter{tocdepth}{1}

\setcounter{tocdepth}{1}
\tableofcontents

\section{Introduction}\label{Intro}

Let $Y$ be a smooth compact connected variety (or orbifold) over an algebraically closed field $\kk$ of characteristic $0$.  The quantum cohomology ring $\QH^*(Y)$ is obtained by using (virtual) counts of holomorphic curves in $Y$ to deform the cup product on the cohomology ring of $Y$.  The degree $0$ subalgebra $\QH^0(Y)$ is very simple, generated by the fundamental class $[Y]=\id\in \QH^*(Y)$.  But now suppose we have a log variety or log orbifold $Y^{\dagger}$ obtained by equipping $Y$ with a log structure, e.g., the data of a reduced effective normal crossings divisor $D\subset Y$.  We call such $(Y,D)$ a \textbf{log pair}.  Then already in codimension $0$ the structure is far more complicated, e.g., logarithmic analogs of the degree $0$ cohomology are often infinite-dimensional.  The question of whether degree $0$ log classes admit an analog of the quantum cohomology product is highly non-trivial, and the resulting algebra $\QH_{\log}^0(Y^{\dagger})$ is expected to be very rich.  Indeed, if $Y^{\dagger}=(Y,D)$ is a \textbf{log Calabi-Yau orbifold with maximal boundary} (i.e., $D$ contains a $0$-stratum and is in $|-K_Y|$), then the Frobenius structure conjecture \cite[arXiv v1, Conj. 0.8]{GHK1} predicts that $\Spec \QH_{\log}^0(Y^{\dagger})$ is the mirror to $Y$.  Furthermore, $\QH_{\log}^0(Y^{\dagger})$ should be naturally equipped with a canonical basis of ``theta functions'' which should agree with the theta functions constructed combinatorially by Gross, Hacking, Keel, Kontsevich, and Siebert \cite{CPS,GHK1,GHKK,GHS}.  Our main result is a proof of this conjecture for cluster varieties, along with a proof for many cases (including all cases with $Y\setminus D$ affine) that the relevant log Gromov-Witten invariants are enumerative.

Cluster varieties were defined in \cite{FG1}, giving geometric meaning to the cluster algebras of \cite{FZ}.  Some examples of cluster varieties include  Grassmannians \cite{ScottGCA} and other partial flag varieties \cite{GLS}, double Bruhat cells of reductive Lie groups \cite{BFZ}, various moduli of local systems (higher-Teichm\"uller spaces) \cite{FG0}, and all two-dimensional log Calabi-Yau varieties with maximal boundary \cite[\S 5]{GHK3}.  By \cite[\S 3.2]{GHK3}, compactifications of cluster varieties can always be viewed as log Calabi-Yau orbifolds with maximal boundary, obtained from toric varieties by preforming certain non-toric blowups of hypertori in the boundary, cf. \S \ref{nontoric}.  \cite{GHKK} constructed canonical theta bases on cluster varieties, and in the process settled many long-standing conjectures about cluster algebras.  We will focus on \textbf{cluster log pairs} $(Y,D)$ as in Definition \ref{clp}.  These are the log pairs obtained by compactifying ``leaves'' of cluster $\s{X}$-varieties, which in our construction includes the usual symplectic leaves and also entire cluster $\s{X}$-spaces, cf. Remark \ref{fullX}.  All theta functions of \cite{GHKK} can be recovered from these cases, cf. Remark \ref{AllTheta}.  This includes the theta functions for log Calabi-Yau surfaces (i.e., Looijenga pairs) as in \cite{GHK1}, cf. Theorem \ref{GHKGHKK}.

\subsection{Relation to other works}\label{other}

Before carefully stating our main results, we note that there are two contemporaneous results \cite{KY,GSInt2} which give other constructions of $\QH^0_{\log}(Y^{\dagger})$ in different generality.  We briefly summarize these works and contrast their results with those of the present paper.

In \cite{KY}, Keel and Yu prove an enumerative version of the Frobenius structure conjecture whenever $Y\setminus D$ is an affine variety and contains a maximal-dimensional algebraic torus.  Their approach is based on Berkovich analytic disks rather than log Gromov-Witten theory, and they obtain descriptions of all the structure constants in terms of these disks.  A priori, their enumerative invariants might be different from the corresponding descendant log Gromov-Witten invariants, but it follows from Theorem \ref{MainNaive} that the invariants agree in all cases where both results apply, i.e., for the subset of cluster varieties which are smooth as varieties (not just as Deligne-Mumford stacks) and which have affine interiors.

As announced in \cite[Thm. 2.2]{GSInt}, Gross and Siebert \cite{GSInt2} give another construction of $\QH_{\log}^0(Y^{\dagger})$ in great generality.  The construction we consider (the Frobenius structure conjecture) involves defining certain $s$-point functions $\QH^0_{\log}(Y^{\dagger})^{\otimes s}\rar \kk\llb \NE(Y)\rrb$, proving that these are given by descendant log Gromov-Witten invariants, and proving that these $s$-point functions for $s=2,3$ uniquely determine the multiplication rule.  This $s$-point function corresponds to taking the constant coefficient---i.e., the coefficient of $\vartheta_0:=1$ in the theta basis expansion---in a product of $s$ functions.   On the other hand, \cite{GSInt2} does not show that the $s$-point functions uniquely determine the multiplication and they do not consider $s>3$.  Instead, they explicitly define all the structure constants in terms of newly developed punctured invariants \cite{ACGSpuncture}, an extension of log invariants which satisfy a powerful splitting lemma---we do not need punctured invariants because we restrict to the $\vartheta_0$ structure constants.  With this approach, Gross and Siebert define an associative algebra associated to log pairs $(Y,D)$ (without the orbifold points that many cluster varieties have) whenever either $K_Y+D$ or $-(K_Y+D)$ is nef, cf. \cite[Thm. 1.9]{GSInt2}.  In particular, this includes the log Calabi-Yau cases $K_Y+D=0$ (in fact, \cite{KY} and \cite{GSInt2} also have results for a weaker version of log Calabi-Yau varieties, cf. \cite[Thm. 1.12]{GSInt2}).   \cite{GSInt2} does not show that their invariants are often enumerative, nor that their theta functions recover those considered in \cite{GHKK}.

We note that our approach is in a sense the reverse of the approach used in \cite{KY,GSInt2}.  In both of those papers, one defines the structure constants directly in terms of curve counting data, and then one proves that the resulting algebra is associative.  For us, the theta functions are already constructed as elements of an associative algebra, and the focus is instead on proving that the $\vartheta_0$ structure constants can be interpreted in terms of curve counts.

We note that \cite{FWY} also seems to be related.  There, one considers the case of smooth $D$ and works with relative invariants, possibly with negative orders of tangency (which also appear in the punctured invariant setting of \cite{GSInt2}).  Fan-Wu-You construct from this a relative version of the quantum cohomology ring.

We briefly give a rough sketch of our arguments---cf. \S \ref{outline} for a more detailed outline.  In \S \ref{Review}, we use \cite[Thm. 3.9]{Man3} to relate the structure constants for the theta functions considered in \cite{GHKK} to certain tropical curve counts, and then we use the main results of \cite{MRud,MRudMult} to relate these tropical invariants to certain descendant log Gromov-Witten invariants of toric varieties.  Here, certain factors in the coefficients are related to monomials in $\kk[\NE(Y)]$ using a new result, Theorem \ref{kappa} (possibly of independent interest), which relates part of the cluster structure data to the lattice of numerical equivalence classes of curves in $Y$.  Then in \S \ref{SectionDegen}, we use a degeneration of $Y^{\dagger}$ and a log Gromov-Witten degeneration formula to relate the invariants of the toric variety to the desired invariants of $Y^{\dagger}$.  In the process, we must show that the relevant log curves are torically transverse, overcoming what was considered to be the main obstacle to generalizing \cite{GPS}.  This involves developing a new tool called ``contractible'' tropical curves, cf. \S \ref{trans}.  These are expansions of ordinary tropical curves which capture the ways that the corresponding log curves might degenerate.  This toric transversality is also essential for showing that the log Gromov-Witten invariants are enumerative in \S \ref{naiveSection}.

\subsection{Statement of the main results}

\subsubsection{Naive counting version}

For a log pair $(Y,D)$, a simple toric blowup $\eta:(\wt{Y},\wt{D})\rar (Y,D)$ is a blowup $\eta:\wt{Y}\rar Y$ of $Y$ along a stratum of $D$, with $\wt{D}$ the reduced inverse image of $D$.  A \textbf{toric blowup} is then a sequence of simple toric blowups.  For $\eta_i:(Y_i,D_i)\rar (Y,D)$, $i=1,2$ two toric blowups of $(Y,D)$, we say that an irreducible component $D'_1\subset D_1$ is equivalent to an irreducible component $D'_2\subset D_2$ if they correspond to the same valuation on the function field of $Y$.   Let $H^0_{\log}(Y,D,\bb{Z})$ denote the free Abelian group generated by $[Y]$ and $[kD']$ for $k\in \bb{Z}_{>0}$ and $D'$ an irreducible component of $\wt{D}$ for some toric blowup $(\wt{Y},\wt{D})$, up to equivalence.  We refer to these generators $[Y]$ and $[kD']$ as \textbf{prime fundamental classes}.\footnote{We note that our prime fundamental classes do indeed form a basis for the degree $0$ classes in the log Chow group $A^0_{\log}(Y^{\dagger})$ as in \cite{Bar,Herr}.}  
 Let $\NE(Y)$ be the cone of effective curve classes in $Y$ up to numerical equivalence (cf. \S \ref{N1}).  Let $\QH_{\log}^0(Y,D)$ denote the completion of $\?{\QH}^0_{\log}(Y,D):=H^0_{\log}(Y,D,\bb{Z})\otimes \kk[\NE(Y)]$ with respect to the unique maximal monomial ideal of $\kk[\NE(Y)]$, so $\QH_{\log}^0(Y,D)$ has the structure of a $\kk \llb \NE(Y)\rrb$-module.  We say that a complete curve $C$ in $Y$ is an \textbf{interior curve} if it is disjoint from $D$, and we say $(Y,D)$ is \textbf{interior-curve free} if it contains no interior curves.  For example, $(Y,D)$ is interior-curve free whenever $Y\setminus D$ is affine.\footnote{In addition to the affine cases, generic two-dimensional log Calabi-Yau varieties with maximal boundary are interior-curve free by \cite[Prop. 4.1]{GHK_MLP}, so Theorem \ref{MainNaive} applies to all the theta functions constructed in \cite{GHK1}.   Interior-curve freeness is not needed for the Gromov-Witten version of our main result, Theorem \ref{MainThm}.}   We say $(Y,D)$ \textbf{supports an ample divisor on its boundary} if there is a toric blowup $(\wt{Y},\wt{D})$ such that $\wt{D}$ supports an effective ample divisor.

\begin{thm}[Main theorem, naive counting version]\label{MainNaive}
 Let $(Y,D)$ be an interior-curve free cluster log pair as in Def. \ref{clp}.  For $\vartheta_1,\ldots,\vartheta_s$ prime fundamental classes of $H^0_{\log}(Y,D,\bb{Z})$, let $\eta:(\wt{Y},\wt{D})\rar (Y,D)$ be a toric blowup in which each $\vartheta_i$ is either $[Y]$ or is represented by $[k_iD_i]$ for $D_i$ an irreducible component of $\wt{D}$.  Let $\beta\in \NE(\wt{Y})$.  Let $(C,x_1,\ldots,x_s,x_{s+1},x_{s+2})$ be a fixed generic irreducible genus $0$ curve with $s+2$ marked points.  Let $y$ be a generically specified point of $\wt{Y}\setminus \wt{D}$.  Define $N_{\beta}^{\naive}(\vartheta_1,\ldots,\vartheta_s)$ to be the number of isomorphism classes of maps $\varphi:C\rar \wt{Y}$ such that $\varphi(x_{s+1})=y$ and $\varphi^* \s{O}_{\wt{Y}}(\wt{D})=\s{O}_C(\sum_{i=1}^s k_ix_i)$, where we take $k_i=0$ if $\vartheta_i=[Y]$.\footnote{This condition means that, for $i=1,\ldots,s$, if $\vartheta_i=[k_iD_i]$, then $\varphi(C)$ intersects $D_i$ at $x_i$ with order $k_i$, and furthermore, these account for all intersections of $\varphi(C)$ with $\wt{D}$.}  All such maps are torically transverse.\footnote{A curve $\varphi:C\rar \wt{Y}$ is torically transverse if $\varphi(C)$ is disjoint from all codimension-two strata of $\wt{D}$}.  Define a $\kk\llb \NE(Y)\rrb$-multilinear $s$-point function $\langle \cdot \rangle^{\naive}:\QH_{\log}^0(Y)^{s} \rar\kk\llb \NE(Y)\rrb$ via
 \begin{align*}
     \langle \vartheta_1,\ldots,\vartheta_s\rangle^{\naive}:=\sum_{\beta\in \NE(\wt{Y})} z^{\eta_*(\beta)}N^{\naive}_{\beta}(\vartheta_1,\ldots,\vartheta_s).
 \end{align*}
Then there is a unique associative product $*$ on $\QH_{\log}^0(Y,D)$ making it into an associative $\kk\llb \NE(Y)\rrb$-algebra such that 
\begin{align*}
    \langle \vartheta_1,\ldots,\vartheta_s\rangle^{\naive} = \langle \vartheta_1*\cdots *\vartheta_s\rangle^{\naive}
\end{align*}
for all $s$-tuples $\vartheta_1, \ldots, \vartheta_s$, $s\geq 1$ (and in fact, the $s=2$ and $s=3$ cases are sufficient to determine $*$).  The algebra $\QH_{\log}^0(Y,D)$ is commutative with identity $[Y]$.  
 If $(Y,D)$ supports an ample divisor on its boundary, then $*$ restricts to give a $\kk[\NE(Y)]$-algebra structure on $\?{\QH}^0_{\log}(Y,D)$. 
 
The algebra $\QH^0_{\log}(Y,D)$ is naturally a subalgebra of $\Gamma(\wt{\s{X}}_{\SSS^{\vee}},\s{O}_{\wt{\s{X}}_{\SSS^{\vee}}})$, where $\wt{\s{X}}_{\SSS^{\vee}}$ denotes a formal version of the Langlands dual $\s{X}$-space\footnote{This is similar to the spaces from the formal Fock-Goncharov conjecture of \cite[\S 6]{GHKK}, but with a different formal completion.} (cf. Remark \ref{Lang}), and the prime fundamental classes are identified with theta functions constructed as in \cite{GHKK}.  If $(Y,D)$ supports an ample divisor on its boundary, then $\?{\QH}^0_{\log}(Y,D)=\Gamma(\?{\s{X}}_{{\SSS}^{\vee}},\s{O}_{\?{\s{X}}_{{\SSS}^{\vee}}})$, where $\?{\s{X}}_{{\SSS}^{\vee}}$ is a partial compactification of the Langlands dual $\s{X}$-space (cf. Remark \ref{Lang}), and again, the prime fundamental classes are the \cite{GHKK} theta functions.
\end{thm}

\subsubsection{The Frobenius structure conjecture}\label{FSC}

For log Calabi-Yau varieties with maximal boundary $(Y,D)$, the Frobenius structure conjecture \cite[arXiv v1, Conj. 0.8]{GHK1} predicts the existence of an algebra structure on $\QH^0_{\log}(Y,D)$, defined essentially as in Theorem \ref{MainNaive} but using descendant log Gromov-Witten invariants in place of the naive curve counts above.  We note that a similar construction from the symplectic perspective was previously suggested by Seidel \cite{Sei1,Sei2}.

Let $(Y,D)$ be a \textbf{smooth} connected log orbifold with maximal boundary.  Here, in addition to $Y$ being smooth as an orbifold (i.e., being a smooth integral separated Deligne-Mumford stack which is proper and finite type over $\kk$), we require that there exists a toric blowup $(Y',D')$ of $(Y,D)$ such that $Y'$ is smooth along $D'$ as a variety (i.e., any orbifold points in the boundary are resolvable by toric blowups), and such that the components of $D'$ are smooth (i.e., $D'$ is snc).  For such a $(Y',D')$, let $S$ denote the dual intersection complex of $D'$.  That is, if $D'=D_1+\ldots+D_n$, then $S$ is the simplicial complex with vertices $v_1,\ldots,v_n$, and with one $(k-1)$-cell with vertices $\{v_{i_1},\ldots v_{i_k}\}$ for each non-empty stratum $D_{i_1}\cap \cdots \cap D_{i_k}$.  Let $B$ be the cone over $S$, and let $\Sigma'$ be the induced simplicial fan in $B$.

In the cone $\sigma\in \Sigma'$ spanned by $v_{i_1},\ldots v_{i_k}$, we have a set of integer points $\sigma(\bb{Z})$ defined as the $\bb{Z}_{\geq 0}$-span of $v_{i_1},\ldots v_{i_k}$ in $\sigma$.  Let $B(\bb{Z}):=\bigcup_{\sigma\in \Sigma'} \sigma(\bb{Z})$ be the integer points of $B$.  Note that there is a bijection between points $p\in B(\bb{Z})$ and prime fundamental classes $\vartheta_p\in \QH^0_{\log}(Y,D)$.  The apex $0\in B(\bb{Z})$ corresponds to $\vartheta_0:=[Y]\in \QH^0_{\log}(Y,\bb{Z})$.  For nonzero $p\in B(\bb{Z})$, we write $\vartheta_p=[|p|D_p]$, i.e., the ray through $p$ determines an irreducible component $D_p$ in the boundary of some toric blowup (up to equivalence), and $|p|\in \bb{Z}_{>0}$ is the \textbf{index} of $p$ in $B(\bb{Z})$---i.e., $|p|$ is the largest positive integer $k$ such that $p=kp'$ for some $p'\in B(\bb{Z})$.  For notational convenience, we write $D_0:=Y$.

By a \textbf{tropical degree}, we mean a map $\Delta:J\rar B(\bb{Z})$ for some finite index-set $J$.  
 In particular, for $\pp$ an $s$-tuple of points $p_1,\ldots,p_s\subset B(\bb{Z})$, we consider 
 \begin{align}\label{Deltapp}
     \Delta_{\pp}:\{1,\ldots,s,s+1,s+2\}\rar B(\bb{Z})
 \end{align} defined by $\Delta_{\pp}(i)=p_i$ for $i=1,\ldots,s$ and $\Delta_{\pp}(s+1)=\Delta_{\pp}(s+2)=0$.  
 
Given any tropical degree $\Delta:J\rar B(\bb{Z})$ (not necessarily $\Delta_{\pp}$), let $\wt{Y}^{\dagger}=(\wt{Y},\wt{D})$ be a toric blowup of $(Y',D')$ with $\wt{Y}$ projective and smooth along $\wt{D}$ as a variety, and such that each $[D_{\Delta(j)}]$ with $\Delta(j)\neq 0$ is represented by an irreducible component of $\wt{D}$.  For $\beta\in \NE(\wt{Y})$, let $\s{M}_{0,\Delta}^{\log}(\wt{Y}^{\dagger},\beta)$ denote \cite{GSlog,AC}'s algebraic moduli stack\footnote{Since any orbifold points of $\wt{Y}$ are away from the boundary, one can use \cite{AGV} to extend the construction of the relevant moduli stacks and their intersection theory to our orbifold setting.  See \cite[\S 5.5]{GPS} for similar considerations from the viewpoint of relative stable maps.} of basic/minimal stable log maps $\varphi^{\dagger}:C^{\dagger}\rar \wt{Y}^{\dagger}$ over $\Spec \kk$ satisfying the following collection of conditions:
\begin{itemize}
\item $C$ has genus $0$;
\item $\varphi_*[C] = \beta$;
\item $C^{\dagger}$ has $|J|$ marked points $\{x_i\}_{i\in J}$ (e.g., in the case of $\Delta=\Delta_{\pp}$, we have $|J|=s+2$);
\item For each $i\in J$, $\varphi(x_i)\in D_{\Delta(i)}$.  Furthermore, if $t_1$ is the generator for the ghost sheaf of $\wt{Y}^{\dagger}$ at a generic point of $D_{\Delta(i)}$, and $t_2$ is the generator for the ghost sheaf of $C^{\dagger}$ at $x_i$, then $\varphi^{\flat}:t_1\mapsto |\Delta(i)|t_2$ where $|\Delta(i)|$ denotes the index of $\Delta(i)$ in $B(\bb{Z})$.  We view this condition as being satisfied automatically for $i$ such that $\Delta(i)=0$.
\end{itemize} 
When the component of $C$ containing $x_i$ is not mapped entirely into $\wt{D}$, this last condition means that the intersection multiplicity of $\varphi(C)$ with $D_{\Delta(i)}$ at $x_i$ is equal to $|p_i|$.  By \cite[Thm. 0.3]{GSlog}, the algebraic stack $\s{M}_{0,\Delta}^{\log}(\wt{Y}^{\dagger},\beta)$ is equipped with a virtual fundamental class $[\s{M}^{\log}_{0,\Delta_{\pp}}(\wt{Y}^{\dagger},\beta)]^{\vir}$ of virtual dimension $$\vdim(\s{M}^{\log}_{0,\Delta_{\pp}}(\wt{Y}^{\dagger},\beta))=\dim(\wt{Y})-(K_Y+D)\cdot \beta+|J|-3,$$
and this class leads to a definition of log Gromov-Witten invariants with the expected properties.

Let $\ev_{i}:\s{M}_{0,\Delta}^{\log}(\wt{Y}^{\dagger},\beta)\rar \wt{Y}$ be the evaluation map $[\varphi^{\dagger}:C^{\dagger}\rar \wt{Y}]\mapsto \varphi(x_i)$.   Let $\pi:\s{C}\rar \s{M}^{\log}_{0,\Delta}(\wt{Y}^{\dagger},\beta)$ denote the universal curve over the moduli space.  Let $\omega_{\pi}$ denote the relative cotangent bundle of $\pi$, and let $\sigma_i$ denote the section of $\pi$ corresponding to $x_i$.  Define
\begin{align}\label{psic1}
    \psi_{i}:=c_1(\sigma_i^* \omega_{\pi}),
\end{align}
i.e., $\psi_i$ is the first Chern class of the line bundle whose fiber over a point $[\varphi^{\dagger}:C^{\dagger}\rar \wt{Y}^{\dagger}]$ is the cotangent space to $C$ at $x_i$.

\begin{dfn}\label{Nbeta}
\begin{align*}
    N_{\beta}(p_1,\ldots,p_s) := \int_{[\s{M}^{\log}_{0,\Delta_{\pp}}(\wt{Y}^{\dagger},\beta)]^{\vir}} \ev_{s+1}^*[\pt] \cdot \psi_{s+1}^{s-1}.
\end{align*}
Note that this log Gromov-Witten invariant will vanish for dimension reasons unless $(K_Y+D)\cdot \beta = 0$.
\end{dfn}

We next define a $\kk\llb \NE(Y)\rrb$-multilinear $s$-point function $$\langle \cdot \rangle:\QH_{\log}^0(Y)^{s} \rar\kk\llb \NE(Y)\rrb$$ via
 \begin{align}\label{spointfunction}
     \langle \vartheta_{p_1},\ldots,\vartheta_{p_s}\rangle:=\sum_{\beta\in \NE(\wt{Y})} z^{\eta_*(\beta)}N_{\beta}(\vartheta_{p_1},\ldots,\vartheta_{p_s})
 \end{align}
where $\eta$ is the composition of blowdowns $\wt{Y}\rar Y'\rar Y$.

\begin{conj}[The Frobenius structure conjecture]\label{FrobConj}
For any smooth connected log Calabi-Yau\footnote{The log Calabi-Yau condition is part of the original Frobenius structure conjecture in \cite[arXiv v1, Conj. 0.8]{GHK1}.  However, as discussed in \S \ref{other}, \cite[Thm. 1.9]{GSInt2} shows that a modified version of the conjecture holds even under the much weaker assumption that $K_Y+D$ is nef or anti-nef.} orbifold with maximal boundary $(Y,D)$, there is a unique associative product $*$ on $\QH_{\log}^0(Y,D)$ making $\QH_{\log}^0(Y,D)$ into an associative $\kk\llb\NE(Y)\rrb$-algebra such that 
\begin{align}\label{spoint}
    \langle \vartheta_{p_1},\ldots,\vartheta_{p_s}\rangle = \langle \vartheta_{p_1}*\cdots *\vartheta_{p_s}\rangle
\end{align}
for all $s$-tuples $p_1, \ldots, p_s\in B(\bb{Z})$, $s\geq 1$.  Furthermore, the $s=2$ and $s=3$ cases of \eqref{spoint} are sufficient to determine $*$.  The algebra $\QH_{\log}^0(Y,D)$ is commutative with identity $\vartheta_0$.  If $(Y,D)$ supports an ample divisor on its boundary, then $*$ restricts to give a $\kk[\NE(Y)]$-algebra structure on $\?{\QH}^0_{\log}(Y,D)$.
\end{conj}

\begin{rmk}\label{logFCA}
The statement of the Frobenius structure conjecture in \cite[arXiv v1, Conj. 0.8]{GHK1} defines $N_{\beta}$ slightly differently.  Their log curves do not include the marked point $x_{s+2}$ (they have only $s+1$ marked points), and their $\psi$-class $\psi_{s+1}$ as in Definition \ref{Nbeta} is only raised to the power of $s-2$, not $s-1$.  That these two definitions are equivalent follows from the Fundamental Class Axiom (generalized to this log setting by the same argument as in the non-log setting).  The advantage of our version of $N_{\beta}$ is that it makes sense for $s=1$ and thus makes the setup more elegant.  We similarly have an extra marked point and $\psi$-class factor in Theorem \ref{QHthm} to allow for the $s=1$ cases.
\end{rmk}

\begin{thm}[Main theorem, Gromov-Witten counting version]\label{MainThm}
Conjecture \ref{FrobConj} holds for all cluster log pairs $(Y,D)$ (as in Def. \ref{clp}).  Furthermore, as in Theorem \ref{MainNaive}, the resulting algebra $\QH^0_{\log}(Y,D)$ is naturally a subalgebra of $\Gamma(\wt{\s{X}}_{\SSS^{\vee}},\s{O}_{\wt{\s{X}}_{\SSS^{\vee}}})$, and the prime fundamental classes are identified with theta functions constructed as in \cite{GHKK}.  If $(Y,D)$ supports an ample divisor on its boundary, then $\?{\QH}^0_{\log}(Y,D)=\Gamma(\?{\s{X}}_{{\SSS}^{\vee}},\s{O}_{\?{\s{X}}_{{\SSS}^{\vee}}})$, where $\?{\s{X}}_{{\SSS}^{\vee}}$ is a partial compactification of the Langlands dual $\s{X}$-space, and again, the prime fundamental classes are the \cite{GHKK} theta functions.
\end{thm}

Theorem \ref{MainNaive} actually follows from Theorem \ref{MainThm} via Proposition \ref{naive}, which says that, under the interior-curve free assumption, the descendant log Gromov-Witten counts $N_{\beta}(p_1,\ldots,p_s)$ agree with the corresponding naive counts $N^{\naive}_{\beta}(\vartheta_{p_1},\ldots,\vartheta_{p_s})$.

As a sample application, recall from \cite{CCGGK} that mirror symmetry for Fano manifolds predicts the equality of the ``quantum period'' of the Fano and the ``classical period'' of a mirror Landau-Ginzburg potential.  This equivalence is a key tool in the ongoing Fano classification program outlined in loc cit.  We will show in a separate paper \cite{ManFano} that this mirror equivalence of quantum and classical periods follows from the Frobenius structure conjecture, at least whenever the log curves being counted are torically transverse.  Our results then imply this equivalence of mirror periods for Fano cluster varieties, thus yielding an algebro-geometric analog of \cite[Thm. 1.1]{Tonk}.

\subsection{Outline of the paper}\label{outline}
In \S \ref{nontoric} we review \cite{GHK3}'s realization of (compactified) cluster varieties as blowups of toric varieties.  Then in \S \ref{N1} we describe the lattice $N_1(Y_{\SSS})$ of curve classes of a cluster variety $Y_{\SSS}$, and in \textbf{Theorem \ref{kappa}} we show how to identify $N_1(Y_{\SSS})$ with the kernel of the exchange matrix of the cluster data.  This result may be of independent interest.

In \S \ref{ScatterIntro} we review the construction of scattering diagrams and theta functions in the context used by \cite{GHKK}.  \textbf{Lemma \ref{nondegen}} (taken from \cite[Thm. 2.17]{Man3}) is essentially the statement that a certain multilinear $s$-point function $\Tr^s$ is sufficient to determine the multiplication rule for the theta functions.  Much of the work of this paper is then to show that the $s$-point function $\Tr^s$ agrees with the $s$-point function $\langle \cdot \rangle$ of \eqref{spointfunction}.

It follows from Theorem \ref{kappa} that the ring of theta functions can be viewed as an algebra over a completion of $\kk[N_1(Y_{\SSS})]$.  In \textbf{Lemma \ref{NECoeff}} we see that the multiplication in fact restricts to give an algebra over $\kk \llb \NE(Y_{\SSS})\rrb$.  This algebra can be identified with a ring of functions on the Langlands dual cluster variety, cf. \textbf{Remark \ref{LangRmk}}.  We take a short detour in \S \ref{GHK1Section} to explain how the constructions of \cite{GHKK} relate to those of \cite{GHK1}---consequently, our main results apply to the \cite{GHK1} theta functions for log Calabi-Yau surfaces (Looijenga pairs) as well.

Next, \S \ref{Review} reviews past results relating theta functions to tropical curve counts (\S \ref{ScatTheta}) and tropical curve counts to descendant log Gromov-Witten invariants of toric varieties (\S \ref{SectionLogGW}).  The combination of these results immediately yields \textbf{Lemma \ref{GWToric}}.  Along the way, we prove (\textbf{Corollary \ref{AmplePolynomial}}) that if $(Y,D)$ supports an ample divisor on its boundary, then the multiplication is polynomial over $\kk[\NE(Y_{\SSS})]$ (i.e., we can restrict to $\?{\QH}^0_{\log}$), thus proving another piece of Theorem \ref{MainThm}.

We prove Theorem \ref{MainThm} in \S \ref{DegenS}.  The idea is to use a degeneration formula to relate the toric descendant log GW invariants from Lemma \ref{GWToric} to the desired descendant log GW invariants of the cluster variety, cf. \textbf{Proposition \ref{DegenProp}}.  The main technical hurdle here is showing that the relevant log stable maps are all torically transverse.  This is achieved in \S \ref{trans} (cf. \textbf{Lemma \ref{DaggerCircLem}}) by proving that the tropicalizations of such log stable maps must be supported on the one-skeleton of the fan for the toric variety.

This toric transversality problem we overcome is perhaps the main reason that \cite{GPS} has not previously been generalized to higher dimensions.  The argument in \S \ref{trans} introduces new ideas, like the concept of a ``contractible'' tropical curve.  Roughly, if $\Gamma$ is a tropicalization of a stable log map $\varphi^{\dagger}$, then the contractible tropical curves which arise as ``expansions'' of $\Gamma$ keep track of the possible tropical types for degenerations of $\varphi^{\dagger}$.  We show that if there were a $\varphi^{\dagger}$ satisfying the generically specified insertions which is not torically transverse, then the tropical type for any deformation of $\varphi^{\dagger}$ satisfying a certain deformation of the insertions would have tropical multiplicity $0$, yielding a contradiction.

In \textbf{Proposition \ref{naive}} we show that $\langle \cdot \rangle$ and $\langle \cdot \rangle^{\naive}$ agree whenever the cluster variety is interior-curve free, thus completing the proof of Theorem \ref{MainNaive}.  The toric transversality result mentioned above is crucial here.

As suggested at the start of the introduction and in the notation, we view $\QH^0_{\log}(Y,D)$ as the degree $0$ part of a conjectural log version of quantum cohomology.  We explain this viewpoint in the appendix by giving a new description of the usual small quantum cohomology ring (\textbf{Theorem \ref{QHthm}}).

\subsection{Acknowledgements}

I am very grateful to Sean Keel and Tony Yue Yu, and to Mark Gross and Bernd Siebert, for keeping me informed about their progress on their closely related projects (cf. \S \ref{other}), as well as for discussions that helped with various technical obstacles which arose in the writing of this and related papers.  In particular, Sean Keel and Mark Gross encouraged the author to pursue this project early on and provided helpful feedback on a draft of this paper.  Collaborations and discussions with Helge Ruddat were also invaluable.  I have also benefited from discussions with Lawrence Barrott, Francesca Carocci, Y.P. Lee, and Dan Pomerleano.  I would also like to thank the anonymous referee for their careful reading and many helpful suggestions.

\section{Cluster varieties}\label{clusterSection}

Here we review the construction of cluster varieties from \cite[\S 1.2]{FG1} as reinterpreted in \cite[\S 3.2]{GHK3}.

\begin{ntn}\label{LatticeNotation}
 For any lattice $L$, let $L_{\bb{Q}}:=L\otimes \bb{Q}$, $L_{\bb{R}}:=L\otimes \bb{R}$, and $T_L:=L\otimes \kk^*$.  We say a nonzero element $v\in L$ is \textbf{primitive} if it is not a positive multiple of any other element of $L$.  Also, for any nonzero $v\in L$, we let $|v|$ denote the \textbf{index} of $v$, i.e., $|v|$ is the unique positive integer such that $v$ is equal to $|v|$ times a primitive vector $v'\in L$.  We use angled brackets $\langle \cdot, \cdot\rangle$ to denote the pairing between a lattice and its dual.
\end{ntn}

\subsection{Construction via blowups of toric varieties}\label{nontoric}

A \textbf{seed} is a collection of data $\SSS$ of the form
\begin{align}\label{seed}
\SSS:=(N,I,E:=\{e_i\}_{i\in I},F\subset I,B),
\end{align}
where $N$ is a lattice of finite rank, $I$ is an index set with $|I| = \rank(N)$, $E$ is a basis for $N$, $F$ is a subset of $I$, and $B$ is a $\bb{Z}$-valued bilinear pairing on $N$.  If $i\in F$, we say $e_i$ is \textbf{frozen}.  Let $I_{\uf}:=I\setminus F$, and let $N_{\uf}$ be the span of $\{e_i\}_{i\in I_{\uf}}$.  The pairing $B$ is required to have a \textbf{skew-symmetrizable} unfrozen part, meaning that there exists a skew-symmetric pairing $\omega$ on $N_{\uf}$ and a collection of positive rational numbers $\{d_i\}_{i\in I_{\uf}}$ such that \begin{align}\label{skewsym}
    B(e_i,e_j)=d_i\omega(e_i,e_j)
\end{align}
for all $i,j\in I_{\uf}$.  The \textbf{Langlands dual} seed $\SSS^{\vee}$ is obtained by replacing $B$ with its negative transpose $-B^T$ while keeping the rest of the seed data the same.

Let $M:=N^*=\Hom(N,\bb{Z})$.  We have two maps $\pi_1, \pi_2:N\rar M$ given by $n\mapsto B(n,\cdot)$ and $n\mapsto B(\cdot,n)$, respectively.   For $i=1,2$, let $K_i:=\ker \pi_i$, and let $N_i:=N/K_i$, which we identify with $\pi_i(N) \subset M$.  We make the following assumptions, although (4) will be relaxed later:\footnote{We note that most these assumptions can always be achieved by modifying only the frozen parts of the data, i.e., without modifying $I_{\uf}$, $N_{\uf}$, or $B|_{N_{\uf}}$.  The only exception is the assumption in 4(b) that $\pi_2(e_i)\neq 0$ for $i\in I_{\uf}$, but vectors failing this assumption can be forgotten from the seed data without affecting the spaces we consider.}
\begin{asss}\label{assume}
~

\begin{enumerate}
    \item $\pi_2(N)$ is saturated in $M$.
    \item $\pi_2(e_i)$ is primitive in $M$ for each $i\in F$.
    \item $\pi_2(e_i)\neq \pi_2(e_j)$ for distinct $i,j\in F$. 
    \item There exists a non-singular complete projective fan $\Sigma$ in $N_2 \otimes \bb{R}$ such that
    \begin{enumerate}
        \item The rays of $\Sigma$ are precisely the rays generated by the vectors $\pi_2(e_i)$ for $i\in F$.
        \item  For each $i\in I_{\uf}$, $\pi_2(e_i)$ is nonzero and is contained in a ray of $\Sigma$, i.e., there is some $j\in F$ such that $\pi_2(e_j)$ points in the same direction as $\pi_2(e_i)$. 
        \item For distinct $i,j\in I_{\uf}$, the rays generated by $\pi_2(e_i)$ and $\pi_2(e_j)$ are either the same or have no cones of $\Sigma$ in common.
        \item Let $\Sigma_i$ denote the set of cones of $\Sigma$ which have codimension $1$ in $N_2\otimes \bb{R}$ and have supports contained in $\pi_1(e_i)^{\perp}$.  Then $\pi_2(e_i)$ is contained in the interior of $\bigcup_{\sigma\in \Sigma_i}\sigma$.
    \end{enumerate}
\end{enumerate}
\end{asss}
Note that Assumption \ref{assume}(1) implies that $N_1$ and $N_2$ are dual to each other via $\langle \pi_1(n_1),\pi_2(n_2)\rangle = B(n_1,n_2)$.  We will write $N_1$ as $\?{M}$ and $N_2$ as $\?{N}$, so $\?{M}=\Hom(\?{N},\bb{Z})$.  Let $r=\rank(\?{N})=\rank(\?{M})$.

Let us fix a fan $\Sigma$ as in Assumption \ref{assume}(4).  Let $\TV_M(\Sigma)$ and $\TV_{\?{N}}(\Sigma)$ denote the toric varieties associated to $\Sigma$ when viewed as a fan in $M$ or $\?{N}$, respectively.  For each $i\in I$, let $D_{\pi_2(e_i)}$ or simply $D_i$ denote the boundary divisor of $\TV_M(\Sigma)$ or $\TV_{\?{N}}(\Sigma)$ corresponding to the ray through $\pi_2(e_i)$ (whether we mean $D_i\subset \TV_M(\Sigma)$ or $D_i\subset \TV_{\?{N}}(\Sigma)$ should always be clear from context).

For each $i\in I_{\uf}$, let $H_i$ denote the scheme-theoretic intersection of $D_{\pi_2(e_i)}\subset \TV_M(\Sigma)$ with the scheme cut out by $(1+z^{e_i})^{|\pi_2(e_i)|}$.  Let $\s{X}_{\SSS,\Sigma}$ denote the scheme obtained by blowing up $\TV_M(\Sigma)$ along $H_i$ for each $i$, let $D^{\s{X}}_{\SSS,\Sigma}$ denote the proper transform of the toric boundary of $\TV_M(\Sigma)$, and let $E^{\s{X}}_i$ denote the exceptional divisor resulting from blowing up $H_i$.  Here, the $H_i$'s may intersect in codimension $2$, so we must choose an order in which to perform the blowups, taking proper transforms of the $H_i$'s and $E^{\s{X}}_i$'s at each step.  Note that centers of blowups associated to distinct $i,j\in I_{\uf}$ are disjoint if $\pi_2(e_i)$ and $\pi_2(e_j)$ are in different rays, but otherwise the choice of ordering may have a codimension-two effect on the resulting space $\s{X}_{\SSS,\Sigma}$, cf. Remark \ref{BlowupRmk}.

Consider the exact sequence
\begin{align*}
	0 \rar K_2 \rar N \stackrel{\pi_2}{\rar} M \stackrel{\lambda}{\rar} K_1^* \rar 0.
\end{align*}
The surjection $\lambda$ induces a map $\lambda:\TV_M(\Sigma)\rar T_{K_1^*}$, and this lifts to a map $\lambda_{\s{X}}:\s{X}_{\SSS,\Sigma}\rar T_{K_1^*}$.  Let $Y_{\SSS,\Sigma}$, or simply $Y_\SSS$, denote a general fiber\footnote{If all of $B$ is skew-symmetrizable, i.e., determined as in \eqref{skewsym} by a skew-symmetric form $\wt{\omega}$ on $N$ and rational numbers $\{d_i\}_{i\in I}$, then the interior of $\s{X}_{\SSS,\Sigma}$ admits a Poisson structure given by $\{z^{e_i},z^{e_j}\}=d_id_j\wt{\omega}(e_i,e_j)z^{e_i+e_j}$, and the interior of $Y_\SSS$ is a symplectic leaf of this Poisson structure.} of $\lambda_{\s{X}}$, and let $D_{\SSS,\Sigma}$ or simply $D_\SSS$ denote the intersection of this fiber with $D^{\s{X}}_{\SSS,\Sigma}$.  Then $(Y_{\SSS},D_{\SSS})$ is a smooth log Calabi-Yau orbifold with maximal boundary.   Let $E_i:=E^{\s{X}}_{i}\cap Y_{\SSS}$.

\begin{rmk}\label{BlowupRmk}
We note that a  fiber $(Y_\SSS,D_\SSS)$ could alternatively be constructed directly by essentially the same construction used to produce $\s{X}_{\SSS,\Sigma}$.  One simply replaces $\TV_M(\Sigma)$ with $\TV_{\?{N}}(\Sigma)$, and replaces each $H_i$ with $\?{H}_i$, defined to be the scheme-theoretic intersection of $D_{\pi_2(e_i)}\subset \TV_{\?{N}}(\Sigma)$ with the scheme cut out by $(a_i+z^{\pi_1(e_i)})^{|\pi_2(e_i)|}$ for some general $a_i\in \kk^*$.  Then $E_i$ is the exceptional divisor associated to blowing up $\?{H}_i$.  Note that for $\pi_2(e_i)$ and $\pi_2(e_j)$ parallel, the corresponding loci $\?{H}_i$ and $\?{H}_j$ might intersect.  But thanks to Assumption \ref{assume}(4)(d), for general fibers of $\lambda_{\s{X}}$, they do not intersect in codimension-two strata of the boundary.  We can therefore apply \cite[Lem. 3.5(1)]{GHK3} to say that the ordering of the blowups only matters up to codimension at most two.  This codimension-two ambiguity will not be important for us.
\end{rmk}

\begin{rmk}\label{fullX}
Since only the unfrozen part of $B$ is required to be skew-symmetrizable, one can always add frozen vectors to obtain a seed $\SSS'$ such that $Y_{\SSS'}$ is a compactification of the full cluster $\s{X}$-variety $\s{X}_{\SSS}$, as opposed to just a compactification of a fiber of $\lambda_{\s{X}}:\s{X}_{\SSS}\rar T_{K_1^*}$. 
\end{rmk}

These pairs $(Y_{\SSS},D_{\SSS})$ are examples of what we call cluster log pairs.  In fact, our definition of cluster log pairs allows for more general boundary:
\begin{dfn}\label{clp}
Let $\SSS=(N,I,E,F,B)$ be a seed $\SSS$ satisfying Assumptions \ref{assume}(1)-(3).  Suppose there exists another seed $\wt{\SSS}=(\wt{N},\wt{I},\wt{E},\wt{F},\wt{B})$ satisfying all of Assumptions \ref{assume}(1)-(4) for some fan $\wt{\Sigma}$, and such the following hold: $N\subset \wt{N}$, $E\subset \wt{E}$,  $I\subset \wt{I}$, $F\subset \wt{F}$ with $I_{\uf}=\wt{I}_{\uf}$, and $B=\wt{B}|_N$.  Let $\Sigma$ be a complete sub cone-complex of $\wt{\Sigma}$ whose rays are precisely the rays generated by the vectors $\pi_2(e_i)$ for $i\in F$.  Furthermore, suppose that the boundary strata of $(Y_{\wt{\SSS},\wt{\Sigma}},D_{\wt{\SSS},\wt{\Sigma}})$ associated to cones of $\wt{\Sigma}\setminus \Sigma$ can be blown down to obtain another smooth log Calabi-Yau orbifold with maximal boundary, which we denote $(Y_{\SSS,\Sigma},D_{\SSS,\Sigma})$.
 A \textbf{cluster log pair} is a pair $(Y_{\SSS,\Sigma},D_{\SSS,\Sigma})$ obtained in this way.
\end{dfn}

In particular, if $\SSS$ and $\Sigma$ do satisfy Assumptions \ref{assume}, then we can take $\wt{\SSS}=\SSS$ and $\wt{\Sigma}=\Sigma$.

\begin{rmk}\label{DMstack} Note that the blowup loci $H_i$ and $\?{H}_i$ are possibly non-reduced, so even if $\Sigma$ is non-singular, the space $Y_{\SSS,\Sigma}$ may still have orbifold singularities.  This is why we work in the generality of Deligne-Mumford stacks.
\end{rmk}

\begin{rmk}\label{AllTheta}
The constructions of \cite{GHKK} involve first constructing theta functions on the cluster variety with principle coefficients $\s{A}^{\prin}$ and then specializing to obtain theta functions on\footnote{We will not further discuss cluster $\s{A}$-varieties, but these spaces (with general coefficients) can be similarly obtained by blowing up partial compactifications of $T_N$ along loci of the form $D_{e_i}\cap Z(a_i+z^{\pi_1(e_i)})$ for $i\in I_{\uf}$, cf. \cite[\S 3.2]{GHK3}.  The map $\pi_2$ lifts to realize these as the universal torsors over the fibers of $\lambda_{\s{X}}$, cf. \cite[\S 4]{GHK3} and \cite{ManCox}.} $\s{A}$ or $\s{X}$, cf. \cite[\S 7.2]{GHKK}.  Theorem \ref{MainThm} applies to determine the theta functions on $\s{X}^{\prin}$ (possibly with some benign modifications to the frozen parts of the seed data), and since $\s{X}^{\prin}$ is isomorphic to $\s{A}^{\prin}$, one can recover all the theta functions of \cite{GHKK}.
\end{rmk}

\subsection{Curve classes}\label{N1}

Letting $b_i$ denote the map blowing up $\?{H}_i$ as in Remark \ref{BlowupRmk}, let $C_i=b_i^{-1}(p)$ for a single generic point $p\in \?{H}_i$, so $C_i$ is a curve contained in $E_i$.  Taking the proper transform of $C_i$ under any remaining blowups, and then the image under blowdowns of extra boundary components, we get a curve $C_i\subset E_i$ in $Y_\SSS$ satisfying 
\begin{align}\label{CiEi}
    [C_i].[E_i]=-\frac{1}{|\pi_2(e_i)|}.
\end{align}

Let $A_*(Y_{\SSS})$ denote the integral Chow lattice of the smooth Deligne-Mumford stack $Y_{\SSS}$, cf. \cite{EG,Kr}.  Consider $N_1(Y_{\SSS})=A_1(Y_{\SSS})$.  We see that $N_1(Y_{\SSS})$ is generated by classes pulled back from $\TV(\wt{\Sigma})$ and then pushed forward from $Y_{\wt{\SSS}}$ to $Y_{\SSS}$, together with the classes $|\pi_2(e_i)|[C_i]$ for $i\in I_{\uf}$.  Let $\NE(Y_{\SSS})$ denote the cone in $N_1(Y_{\SSS})$ generated by classes of effective curves.

The following useful theorem\footnote{A version of Theorem \ref{kappa} in dimension $2$ without frozen vectors or the corresponding boundary divisors was proven in \cite[Thm. 5.5]{GHK3}.  The argument here is inspired by the proof in loc. cit.} identifies the lattice $K_2\subset N$ with the lattice $N_1(Y_{\SSS,\Sigma})$. 

\begin{thm}\label{kappa} 
There is a unique isomorphism 
\begin{align*}
    \kappa:K_2\risom N_1(Y_{\SSS,\Sigma})
\end{align*}
taking a vector $k=\sum_{i\in I} a_ie_i$ to a curve class $[C]_k$ such that $[C]_k.[E_i]=a_i$ for each $i\in I_{\uf}$, and $[C]_k.[D_i]=a_i$ for each $i\in F$.
\end{thm}
\begin{proof}
Let us first assume that $\SSS=\wt{\SSS}$ and $\Sigma=\wt{\Sigma}$ in the construction of $(Y_{\SSS,\Sigma},D_{\SSS,\Sigma})$ as in Definition \ref{clp}.  Let $$\eta_{\SSS}:Y_{\SSS}\rar \TV(\Sigma):=\TV_{\?{N}}(\Sigma)$$ denote the blowdown map.  By standard toric geometry (cf. \cite[\S 3.4]{Fult}), $A_{n-1}(\TV(\Sigma))$ is spanned by the classes of its boundary divisors. Hence, $A_{n-1}(Y_{\SSS,\Sigma})$ is spanned by the boundary divisor classes $[D_i]$ with $i\in F$, together with the exceptional divisor classes $[E_i]$ with $i\in I_{\uf}$.  So a class $[C]\in N_1(Y_{\SSS,\Sigma})$ is indeed uniquely determined by its intersections with the classes $[D_i]$, $i\in F$ and $[E_i]$, $i\in I_{\uf}$.

We now check that such a $[C]_k$ exists for each $k\in K_2$.  By definition, $\sum_{i\in I} a_ie_i\in K_2$ gives a relation $\sum_{i\in I} a_i\pi_2(e_i) = 0$, and such a relation corresponds to a class $[\?{C}]_k\in N_1(\TV(\Sigma))$ which, for each ray $\rho \in \Sigma$, satisfies 
 \begin{align}\label{Ckbar}
 [\?{C}]_k.[D_{\rho}]=\sum_{\{i\in I:\pi_2(e_i)\in \rho\}} a_i|\pi_2(e_i)|.    
 \end{align}
Now let 
\begin{align*}
[C]_k:=\eta_{\SSS}^*[{\?{C}}]_k-\sum_{i\in I_{\uf}} a_i|\pi_2(e_i)|[C_i].     
\end{align*}
Using \eqref{CiEi}, it is straightforward to check that $[C]_k$ has the desired intersection multiplicities with $[E_i]$ or $[D_i]$ for each $i\in I$.

Next, we want to check that $\kappa$ is surjective.  It is clear that the image includes the pullback of any class from $N_1(\TV(\Sigma))$, so we just have to check that the image includes the classes $|\pi_2(e_i)|[C_i]$ for each $i\in I_{\uf}$.    By Assumption \ref{assume}(2) and (4)(b), there is some $j\in F$ and $a\in \bb{Z}_{>0}$ such that $\pi_2(e_i)=a \pi_2(e_j)$.  For $k=-e_i+ae_j$, the class $[\?{C}]_k$ is just $0\in N_1(\TV(\Sigma))$, so we have
\begin{align*}
 \kappa(-e_i+ a e_j)=|\pi_2(e_i)|[C_i],
\end{align*}
as desired.

We have thus proved the claim when $\SSS=\wt{\SSS}$, $\Sigma=\wt{\Sigma}$.  For the more general situation, let $\wt{\kappa}$ and $\wt{K_2}$ denote the appropriate data associated to $\wt{\SSS}$ and $\wt{\Sigma}$.  Consider the blowdown $$\eta_{\Sigma}:Y_{\wt{\SSS},\wt{\Sigma}}\rar Y_{\SSS,\Sigma}.$$  Note that the inclusion $N\subset \wt{N}$ identifies $K_2$ with $\wt{K}_2\cap N$.  It follows from the projection formula for Chow rings that $\eta_{\Sigma}^* N_1(Y_{\SSS,\Sigma})$ is a sublattice of $N_1(Y_{\wt{\SSS},\wt{\Sigma}})$ and consists of those classes which have $0$ intersection with the boundary divisors $D_i$ for each $i\in \wt{F}\setminus F$.  Hence, $$(\wt{\kappa})^{-1}(\eta_{\Sigma}^* N_1(Y_{\SSS,\Sigma})) = K_2.$$  
We now see (using the projection formula again) that the desired map $\kappa$ is 
\begin{align*}
    \kappa:=(\eta_{\Sigma})_* \circ \wt{\kappa}|_{K_2}.
\end{align*}
\end{proof}

Note that for each $i\in I_{\uf}$, the class $|\pi_2(e_i)|[C_i]\in N_1(Y_{\SSS})$ generates an extremal ray of the Mori cone $\NE(Y_{\SSS})$.  Let $\NE(Y_\SSS)_\SSS$ denote the localization of $\NE(Y_\SSS)$ obtained by adjoining $-|\pi_2(e_i)|[C_i]$ for each $i\in I_{\uf}$.  Let $N^{\oplus}$ denote the submonoid of $N$ spanned by the elements $e_i$, $i\in I$, i.e., 
\begin{align}\label{Noplus}
    N^{\oplus} = \left\{\sum_{i\in I} a_ie_i| a_i\in \bb{Z}_{\geq 0} \mbox{ for each } i\in I\right\}.
\end{align}
Let $K_2^{\oplus}:=K_2\cap N^{\oplus}$.
\begin{lem}\label{NEloc}
$\kappa(K_2^{\oplus})\subset \NE(Y_\SSS)_\SSS$.
\end{lem}
\begin{proof}
This follows from the proof of Theorem \ref{kappa}, noting that when each $a_i\geq 0$, the class $[\?{C}]_k$ of \eqref{Ckbar} can be represented by an effective curve in $\TV(\Sigma)$.
\end{proof}

\begin{rmk}\label{Lang}
The \textbf{Langlands dual $\s{X}$-space} is the space $\s{X}_{\SSS^{\vee}}$ associated to the Langlands dual seed $\SSS^{\vee}$ (with the boundary removed).  It comes with a map $\lambda_{\s{X}}:\s{X}_{\SSS^{\vee}}\rar T_{K_2^*}=\Spec \kk[K_2]$, which by Theorem \ref{kappa} can be viewed as
\begin{align*}
    \lambda_{\s{X}}:\s{X}_{\SSS^{\vee}} \rar \Spec \kk[N_1(Y_{\SSS,\Sigma})].
\end{align*}
By adding appropriate boundary strata, this can be extended to a family 
\begin{align*}
    \lambda_{\s{X}}:\?{\s{X}}_{\SSS^{\vee}} \rar \Spec \kk[\NE(Y_{\SSS,\Sigma})],
\end{align*}
and then one can define the formal completion $\wt{\s{X}}_{\SSS^{\vee}}$ of $\?{\s{X}}_{\SSS^{\vee}}$ at the boundary $\?{\s{X}}_{\SSS^{\vee}}\setminus \s{X}_{\SSS^{\vee}}$.
\end{rmk}

\section{Scattering diagrams and theta functions}

\subsection{Review of scattering diagrams and theta functions}\label{ScatterIntro}  We next recall the notion of a scattering diagram and the construction of theta functions.  We continue to assume we have the data of a seed $\SSS$ and fan $\Sigma$ as in \S \ref{nontoric}.

Recall $N^{\oplus}\subset N$ as in \eqref{Noplus}, and let $N^+:=N^{\oplus}\setminus \{0\}$.  Consider $\kk[N^{\oplus}]$.  Let $\f{m}$ denote the unique maximal monomial ideal of $\kk[N^{\oplus}]$, that is, the ideal generated by all $z^n$ with $n\in N^+$.  For each $k\in \bb{Z}_{>0}$, we can take the quotient $\kk[N^{\oplus}]/\f{m}^k$, and we thus define the inverse limit $\kk\llb N^{\oplus}\rrb:=\varprojlim_k \kk[N^{\oplus}]/\f{m}^k$.

Now let 
\begin{align}\label{P}
    P:=N^{\oplus}+\kappa^{-1}(\NE(Y_{\SSS}))\subset N.
\end{align}  Let $A:=\kk[P]$ with its obvious $P$-grading.  Let $\wh{A}$ denote the $N^+$-adic completion of $\kk[P]$.  I.e., for each $k\in \bb{Z}_{\geq 1}$, denote
\begin{align*}
    kN^+:=\{n_1+\ldots,n_k\in P|n_i\in N^+ \mbox{ for each $i=1,\ldots,k$}\}.
\end{align*}
Then $\wh{A}$ consists of Laurent series $\sum_{p\in P} a_p z^p$ such that, for each $k\in \bb{Z}_{\geq 1}$, $a_p=0$ for all but finitely many $p\in P\setminus kN^+$.  Equivalently, $\wh{A}=\kk\llb N^{\oplus}\rrb \otimes_{\kk[N^{\oplus}]} \kk[P]$, or if we define $A_k:=\kk[N^{\oplus}]/\f{m}^k \otimes_{\kk[N^{\oplus}]} \kk[P]$, then  $\wh{A}=\varprojlim_k A_k$.

Similarly, let $K_2^+:=K_2^{\oplus}\setminus \{0\}$, where we recall $K_2^{\oplus}=K_2\cap N^{\oplus}$.  Let $R:=\kk[K_2\cap P]$, and let $\wh{R}=\kk\llb K_2^{\oplus}\rrb \otimes_{\kk[K_2^{\oplus}]} R$ be the $K_2^+$-adic completion of $R$.  We sometimes view $A$ and $\wh{A}$ as algebras over $R$ and $\wh{R}$, respectively.

Recall our notation $\?{N}:=N_2\subset M$ and $\?{M}:=N_1=\?{N}^*$.

\begin{dfn}
A \textbf{wall} $(\f{d},f)$ in $\?{N}_{\bb{R}}$ is the data of a function $$f\in \kk\llb z^n\rrb\subset \kk\llb N^{\oplus}\rrb$$ for some $n\in N^+$, and a convex (but not necessarily strongly convex) codimension one rational polyhedral cone $$\f{d}\subset \?{N}_{\bb{R}}$$  such that $f\equiv 1$ modulo $\f{m}$ and such that the linear span of $\f{d}$ contains $\pi_2(n)$.  The vector $-\pi_2(n)\in \?{N}$ is called the \textbf{direction} of the wall.  The wall is called \textbf{incoming} if $\f{d}$ contains $\pi_2(n)$ and \textbf{outgoing} otherwise.

A \textbf{scattering diagram} $\f{D}$ is a set of walls in $\?{N}_{\bb{R}}$ such that for each $k >0$, there are only finitely many walls $(\f{d},f)$ with $f$ not equivalent to $1$ modulo $\f{m}^k$.  Given $\f{D}$, we let $\f{D}^k$ denote the finite scattering diagram consisting of walls $(\f{d},f) \in \f{D}$ for which $f\not\equiv 1$ modulo $\f{m}^k$.
\end{dfn}

  Denote $\Supp(\f{D}):= \bigcup_{(\f{d},f)\in \f{D}} \f{d}$, and \begin{align*}
\Joints(\f{D}):= \bigcup_{(\f{d},f)\in \f{D}} \partial \f{d} \cup \bigcup_{\substack{(\f{d}_1,f_1),(\f{d}_2,f_2)\in \f{D}\\
		               \codim (\f{d}_1\cap \f{d}_2\subset \?{N}_{\bb{R}}) = 2}} \f{d}_1\cap \f{d}_2.
\end{align*}
We will sometimes denote a wall $(\f{d},f)$ by just $\f{d}$.  On the other hand, we may write $(\f{d},f\in \kk\llb z^n\rrb)$ if we want to explicitly indicate the data of $n$.

Consider a smooth immersion $\gamma:[0,1]\rar \?{N}_{\bb{R}}\setminus \Joints(\f{D})$ with endpoints not in $\Supp(\f{D})$ which is transverse to each wall of $\f{D}$ it crosses.  Let $(\f{d}_i,f_i\in \kk\llb z^{n_i}\rrb)$, $i=1,\ldots, s$, denote the walls of $\f{D}^{k}$ crossed by $\gamma$, and say they are crossed at times $0<t_1\leq \ldots \leq t_s<1$, respectively (the ambiguity in the labelling when $t_i=t_{i+1}$ is unimportant).   Define a ring automorphism $\theta_{\f{d}_i}$ of $\kk[N^{\oplus}]/\f{m}^k$ which, for each $n\in N^{\oplus}$, acts via
\begin{align}\label{WallCross}
\theta_{\f{d}_i}(z^p):=z^pf_i^{\langle u_i, \pi_2(p)\rangle},
\end{align}
where $u_i$ is the primitive element of $\f{d}_i^{\perp}\subset \?{M}$ which is positive on $-\gamma'(t_i)$. Let $\theta_{\gamma,\f{D}}^k:=\theta_{\f{d}_s} \circ \cdots \circ \theta_{\f{d}_1}\in \Aut(A_k)$. Finally, define the path-ordered product
\begin{align*}
\theta_{\gamma,\f{D}}:= \varprojlim_k \theta_{\gamma,\f{D}}^k \in \Aut(\wh{A}).
\end{align*}

We say two scattering diagrams $\f{D}$ and $\f{D}'$ are  \textbf{equivalent} if $\theta_{\gamma,\f{D}} = \theta_{\gamma,\f{D}'}$ for each smooth immersion $\gamma$ as above.  One says $\f{D}$ is \textbf{consistent} if each $\theta_{\gamma,\f{D}}$ depends only on the endpoints of $\gamma$.

The following theorem of Gross-Siebert and Kontsevich-Soibelman is fundamental to the theory of scattering diagrams.
\begin{thm}[\cite{GS11}, \cite{WCS}]\label{KSGS}
Let $\f{D}_{\In}$ be a finite scattering diagram in $\?{N}_{\bb{R}}$ whose only walls have full hyperplanes as their supports.  Then there is a unique-up-to-equivalence scattering diagram $\f{D}$, also denoted $\scat(\f{D}_{\In})$, such that $\f{D}$ is consistent, $\f{D} \supset \f{D}_{\In}$, and $\f{D}\setminus \f{D}_{\In}$ consists only of outgoing walls.
\end{thm}

Let us now fix a consistent scattering diagram $\f{D}$ in $\?{N}_{\bb{R}}$.  Let $\varphi:\?{N}_{\bb{R}}\rar N_{\bb{R}}$ denote the integral $\Sigma$-piecewise-linear section of $\pi_2$ determined by setting 
\begin{align}\label{varphi}
    \varphi(\pi_2(e_i)):=e_i
\end{align} for each $i\in F$ and then extending linearly over the cones of $\Sigma$.

\begin{dfn}\label{broken line}
Let $p \in \?{N}$, $Q\in \?{N}_{\bb{R}}\setminus \Supp(\f{D})$.  A  broken line $\gamma$ with  ends $(p,Q)$ is the data of a continuous map $\gamma:(-\infty,0]\rar \?{N}_{\bb{R}}\setminus \Joints(\f{D})$, values $t_0 \leq t_1 \leq \ldots \leq t_{\ell-1} < t_{\ell} = 0$, and for each $i=0,\ldots,\ell$, an associated element $c_iz^{v_i}\in \wh{A}$, such that:
\begin{itemize}[noitemsep]
\item $\gamma(0)=Q$.
\item For $i=1\ldots, \ell$, $\gamma'(t)=-\pi_2(v_i)$ for all $t\in (t_{i-1},t_{i})$.  Similarly, $\gamma'(t)=-\pi_2(v_0)$ for all $t\in (-\infty,t_0)$. 
\item $c_0=1$ and $v_0=\varphi(p)$.
\item For $i=0,\ldots,\ell-1$, $\gamma(t_i)\in \f{d}_i$ for some wall $(\f{d}_i,f_i)\in \f{D}$, and $c_{i+1}z^{v_{i+1}}\neq c_iz^{v_i}$ is a monomial term in the power series expansion of $c_iz^{v_i}f_i^{\langle u_i,v_i\rangle}$, where $u_i$ is the primitive element of $\f{d}_i^{\perp}$ which is positive on $v_i$ (i.e., $c_iz^{v_i}f_i^{\langle u_i,v_i\rangle}$ is $\theta_{\f{d}_i}(c_iz^{v_i})$ as defined in \eqref{WallCross} for a smoothing of $\gamma$).
\end{itemize}
\end{dfn}

Fix a generic point $Q \in \?{N}_{\bb{R}}\setminus \Supp(\f{D})$. For any $p\in \?{N}$, we define a theta function
\begin{align}\label{vartheta-dfn}
\vartheta_{p,Q}:=\sum_{\Ends(\gamma)=(p,Q)} c_{\gamma}z^{n_{\gamma}} \in \wh{A}.
\end{align}
Here, the sum is over all broken lines $\gamma$ with ends $(p,Q)$, and $c_{\gamma}z^{n_{\gamma}}$ denotes the monomial attached to the final straight segment of $\gamma$.  In particular, we define $\vartheta_{0,Q}:=1$.  One can prove that these functions $\vartheta_{p,Q}$ form a well-defined topological $\wh{R}$-module basis for $\wh{A}$, hence also for the $\wh{R}$-subalgebra $\wh{A}_{\Theta,Q}\subset \wh{A}$ which they generate, cf. \cite[Prop 2.14]{Man3}.

Furthermore, if $Q$ and $Q'$ are two generic points in $\?{N}_{\bb{R}}\setminus \Supp(\f{D})$, and if $\gamma$ is a smooth path from $Q$ to $Q'$ avoiding $\Joints(\f{D})$, then an important result of \cite{CPS} says that $\vartheta_{p,Q'}=\theta_{\gamma,\f{D}}(\vartheta_{p,Q})$.  Hence, as an abstract algebra, $\wh{A}_{\Theta,Q}$ is independent of $Q$, and so we denote it by just $\wh{A}_{\Theta}$.  Similarly, we denote $\vartheta_{p,Q}\in \wh{A}_{\Theta,Q}=\wh{A}_{\Theta}$ by simply $\vartheta_p$.

Given $f=\sum_{p\in \?{N}} c_p\vartheta_p\in \wh{A}_{\Theta}$, define $\Tr(f):=c_0\in \wh{R}$.  This determines a symmetric multilinear $s$-point function $\Tr^s:\wh{A}_{\Theta}^s \rar \wh{R}$,
\begin{align}\label{Trs}
    \Tr^s(f_1,\ldots,f_s):=\Tr(f_1\cdots f_s).
\end{align}
The following lemma says that these $\Tr^s$ are sufficient to determine the entire multiplication structure.
\begin{lem}[\cite{Man3}, Thm. 2.17]\label{nondegen}
 $\Tr^2$ is non-degenerate as a symmetric $\wh{R}$-bilinear pairing on the $\wh{R}$-module $\wh{A}_{\Theta}$.  Thus, given $\wh{A}_{\Theta}$ as a module over $\wh{R}$ topologically generated by $\{\vartheta_p\}_{p\in \?{N}}$, the multiplication rule giving the $\wh{R}$-algebra structure is uniquely determined by $\Tr^2$ and $\Tr^3$. 
\end{lem}
We will prove, for a certain $\f{D}_{\In}$, that $\Tr^s$ as defined here is given by the $s$-point function of \eqref{spoint}.  Specifically, the initial scattering diagram $\f{D}^\SSS_{\In}$ with which we shall work is
\begin{align}\label{Din}
    \f{D}^\SSS_{\In}:=\left\{\left(\pi_1(e_i)^{\perp},(1+z^{e_i})^{|\pi_1(e_i)|}\right) : i\in I_{\uf}\right\},
\end{align}
and we denote
\begin{align*}
    \f{D}^\SSS:=\scat(\f{D}^\SSS_{\In})
\end{align*}

\subsection{Effectiveness of curve classes}\label{EffCurve}

We have seen that the theta functions form a basis for $\wh{A}_{\Theta}$ over  $\wh{R}$, but Theorems \ref{MainNaive} and \ref{MainThm} claim that we should be able to work over $\kk\llb\NE(Y_{\SSS})\rrb$.  This is the content of the following Lemma:

\begin{lem}\label{NECoeff}
The theta functions $\{\vartheta_p\}_{p\in \?{N}}$ form a topological $\kk\llb \NE(Y_{\SSS})\rrb$-module basis for a $\kk\llb \NE(Y_{\SSS})\rrb$-algebra $A_{\Theta}$.
\end{lem}

\begin{proof}

From the statement that the theta functions form a topological $\wh{R}$-module basis for $\wh{A}_{\Theta}$, we know that for any $p,q\in \?{N}$, we can write
\begin{align*}
    \vartheta_p\vartheta_q=\sum_r c_{pqr} \vartheta_r
\end{align*}
for some collection of ``structure constants'' $c_{pqr}\in \wh{R}$.  Recall that $\wh{R}=\kk\llb K_2^{\oplus}\rrb \otimes_{\kk[K_2^{\oplus}]} \kk[K_2\cap P]$.  By Lemma \ref{NEloc}, $\kappa(K_2^{\oplus}) \subset \NE(Y_{\SSS})_{\SSS}$, hence $\kappa(K_2\cap P) \subset \NE(Y_{\SSS})_{\SSS}$ as well.  Hence, each $c_{pqr}$ is a formal sum of monomials of the form $c_{pqr\beta}z^{\beta}$ with $c_{pqr\beta}\in \kk$ and $\beta\in \NE(Y_{\SSS})_{\SSS}$ (under the identification $\kappa$).

We would like to show that the localization of the Mori cone here is in fact not necessary, i.e., we want to show that $c_{pqr\beta}=0$ unless $\beta\in \NE(Y_{\SSS})\subset \NE(Y_{\SSS})_{\SSS}$.  It then follows that the $\kk\llb \NE(Y_{\SSS})\rrb$-submodule $A_{\Theta}$ of $\wh{A}_{\Theta}$ spanned topologically by the theta functions is in fact a $\kk\llb \NE(Y_{\SSS})\rrb$-subalgebra.

For this, we take advantage of the operation of mutation.  Namely, for each element $e_i\in E$ with $i\in I_{\uf}$ from our seed $\SSS$, we consider the birational map $\mu_i:T_{\?{N}} \dashrightarrow T_{\?{N}}$, $z^{\pi_1(n)}\mapsto z^{\pi_1(n)}(1+z^{\pi_1(e_i)})^{B(e_i,n)}$.   \cite{GHK3} shows that this map can be interpreted geometrically (up to codimension $2$) as taking the blowup of $\?{H}_i$ as in Remark \ref{BlowupRmk}, followed by taking the blowdown of a certain locus $\wt{F}_i\rar \?{H}_i'$ as in Figure \ref{mutation-fig}.  Let $C_i'$ denote a generic fiber of the blowdown $\wt{F}_i\rar \?{H}_i'$, so $|\pi_2(e_i)|[C_i']$ generates another extremal ray of $\NE(Y_{\SSS,\Sigma})$.  Let $\NE(Y_{\SSS})_{\mu_i(\SSS)}$ denote the localization of $\NE(Y_{\SSS})$ obtained by adjoining $-|\pi_2(e_j)|[C_j]$ for $j\in I_{\uf}\setminus \{i\}$, along with $-|\pi_2(e_i)|[C'_i]$.

\begin{figure}[htb]
\center{\resizebox{4.8 in}{1.2 in}{
\xy
(-60,-30)*{}; (-25,-30)*{} **\dir{-}="Dun1";
(-60,-10)*{}; (-25,-10)*{} **\dir{-}="Dup1";
(-40,-8)*{}; (-40,-32)*{} **\dir{-}="F1";
(-60,-27)*{D_{\pi_2(e_i)}};
(-60,-7)*{D_{-\pi_2(e_i)}};
(-42,-20)*{F_i};
(-5,10)*{}; (30,10)*{} **\dir{-}="Dun2";
(-5,30)*{}; (30,30)*{} **\dir{-}="Dup2";
(15,31)*{}; (23,18)*{} **\dir{-}="F2";
(15,9)*{}; (23,22)*{} **\dir{-}="E2";
(-5,13)*{D_{\pi_2(e_i)}};
(-5,33)*{D_{-\pi_2(e_i)}};
(22,26)*{\wt{F}_i};
(22,14)*{\wt{E}_i};
(50,-30)*{}; (85,-30)*{} **\dir{-}="Dun3";
(50,-15)*{}; (85,0)*{} **\dir{-}="Dup3";
(70,-5)*{}; (70,-31)*{} **\dir{-}="E3";
(50,-28)*{D_{\pi_2(e_i)}};
(50,-10)*{D_{-\pi_2(e_i)}};
(68,-20)*{E_i};
{\ar (-11,4);(-19,-4)};
{\ar (36,4);(44,-4)};
(-40.1,-30)*{\bullet};
(-35.7,-32.7)*{\?{H}_i};
(70,-6.5)*{\bullet};
(67,-4)*{\?{H}_i'};
{\ar@{-->}  (-9,-22);(33,-22)};
(12.5,-19.5)*{\mu_i}
\endxy
}}
\caption{Consider the fan $\Sigma_i$ in $\?{N}_{\bb{R}}$ with rays generated by $\pi_2(e_i)$ and $-\pi_2(e_i)$.  Let $\pi_2(e_i)'\in \?{N}$ denote the primitive element with direction $\pi_2(e_i)$.  The map $\?{N}\rar \?{N}/\bb{Z}\pi_2(e_i)'$ induces a $\bb{P}^1$ fibration of the toric variety $\TV(\Sigma_i)$ over $T_{\?{N}/\bb{Z}\pi_2(e_i)'}$. The mutation $\mu_{i}$ is the birational map $T_{\?{N}}\dashrightarrow T_{\?{N}}$ given by including $T_{\?{N}}$ into $\TV(\Sigma_i)$, blowing up the locus $\?{H}_i$ (left arrow),  contracting the proper transform $\wt{F}_i$ of the fibers $F_i$ which hit $\?{H}_i$ down to a locus $\?{H}_i'$ in $D_{-\pi_2(e_i)}$ (right arrow), and then taking the complement of the proper transforms of the boundary divisors.}\label{mutation-fig}
\end{figure} 

There is a well-known seed-mutation which associates a new seed $\SSS_i$ to each $i\in I$ above.  We can use this to define a different scattering diagram $\f{D}^{\mu_i(\SSS)}$ and an associated algebra of theta functions.  It follows from \cite[Thm 1.24]{GHKK} that the resulting theta functions have the same multiplication rule as before.  Thus, for a coefficient $c_{pqr\beta}$ of theta the function multiplication to be nonzero, we must have \begin{align}\label{sNE}
\beta\in \NE(Y_{\SSS})_\SSS\cap \bigcap_{i\in I_{\uf}} \NE(Y_{\SSS})_{\mu_i(\SSS)}.
\end{align}
Since the classes $|\pi_2(e_i)|[C_i]$ and $|\pi_2(e_i)|[C'_i]$ for $i\in I_{\uf}$ are all distinct and extremal in $\NE(Y_{\SSS})$, the intersection \eqref{sNE} is just $\NE(Y_{\SSS})$, as desired.
\end{proof}

\begin{rmk}[Relation to the Langlands dual spaces]\label{LangRmk}
By \cite[\S 3.2]{GHK3}, $\s{X}_{\SSS^{\vee}}$ is covered up to codimension $2$ by the initial cluster and the mutation-adjacent clusters.  The same then follows for $\?{\s{X}}_{\SSS^{\vee}}$ and $\wt{\s{X}}_{\SSS^{\vee}}$.  The birational automorphisms gluing these clusters in $\wt{\s{X}}_{\SSS^{\vee}}$ are precisely the automorphisms associated to crossing the walls of $\f{D}^{\SSS}_{\In}$.  Thus, the theta functions are realized as global regular functions on $\wt{\s{X}}_{\SSS^{\vee}}$.

Now suppose that the theta functions in fact generate a polynomial algebra over $\kk[\NE(Y_{\SSS})]$.  The identification of this algebra with $\Gamma(\?{\s{X}}_{{\SSS}^{\vee}},\s{O}_{\?{\s{X}}_{{\SSS}^{\vee}}})$ then follows from \cite[Thm. 0.3]{GHKK}.
\end{rmk}

\subsection{Relation to \cite{GHK1}}\label{GHK1Section}
We briefly explain how the constructions considered here, i.e., those of \cite{GHKK}, relate to the constructions of \cite{GHK1}.  This subsection will not be used elsewhere in the present paper.

Suppose that $(Y_{\SSS,\Sigma},D_{\SSS,\Sigma})$ is a smooth log Calabi-Yau surface, or in the language of \cite{GHK1}, a Looijenga pair.  By \cite[Prop. 1.3]{GHK1}, every Looijenga pair has a toric blowup which admits a toric model.  Hence, every Looijenga pair can be obtained as $(Y_{\SSS,\Sigma},D_{\SSS,\Sigma})$ for some $\SSS$ and $\Sigma$, as observed in \cite[\S 5]{GHK3}.

\begin{thm}\label{GHKGHKK}
The theta functions constructed from $(Y_{\SSS,\Sigma},D_{\SSS,\Sigma})$ using the procedures of \cite{GHK1} agree with the theta functions constructed as in \S \ref{ScatterIntro}.  I.e., the construction of \cite{GHK1} can be viewed as a special case of the constructions of \cite{GHKK}.
\end{thm}

The key observation behind our proof is that the map $\kappa$ of \eqref{kappa} identifies $\pi_2:N\rar \?{N}$ with the integral points of the local system $\s{P}\rar B$ considered in \cite{GHK1} (after moving worms), and this identifies the piecewise-linear section $\varphi$ of $\pi_2$ as in \eqref{varphi} with the analogous piecewise-linear section in \cite{GHK1}.

\begin{proof}
Since the theta functions behave well with respect to toric blowups, we may assume we have a toric model $\eta_{\SSS}:Y_{\SSS,\Sigma}\rar \TV(\Sigma)$.  The main construction in \cite{GHK1} uses a scattering diagram which lives in a singular integral linear manifold $B$, with coefficients coming from a local system $\s{P}$ over $B$.  But by the ``moving worms'' construction of \cite[\S 3.2-3.4]{GHK1}, we can obtain the same theta functions (restricted to some Zariski dense open subset of $\Spec \kk\llb \NE(Y_{\SSS,\Sigma})\rrb$, i.e., after some localization of the coefficient ring) by using a scattering diagram $\?{\f{D}}$ in $\?{N}_{\bb{R}}$.

Specifically, as in \cite[\S 3.4, the lead up to Theorem 3.33]{GHK1}, $\?{\f{D}}:=\scat(\?{\f{D}}_0)$ for $\?{\f{D}}_0$ defined as follows.  Let $\?{Y}:=\TV(\Sigma)$, $\?{P}:=\NE(\?{Y})\subset A_1(\?{Y},\bb{Z})$, and let $\?{\varphi}:\?{N}_{\bb{R}}\rar A_1(\?{Y},\bb{R})$ be the $\Sigma$-piecewise linear integral strictly $P$-convex function function whose bending parameter along the ray $\rho$ corresponding to $D_{\rho}$ is $[D_{\rho}]$.  Let $\varphi':=\eta_{\SSS}^*\circ \?{\varphi}:\?{N}_{\bb{R}}\rar A_1(Y,\bb{R})$.  For $\varphi$ as in \eqref{varphi} and $\kappa$ as in Theorem \ref{kappa}, we claim that $\varphi(n)=(n,\kappa^{-1}\circ\varphi'(n))$.

To see this, for each $i\in F$, let $\?{e}_i:=\pi_2(e_i)$ and let $\f{b}_i\in A_1(Y,\bb{Z})$ denote the bending parameter of $\varphi$ along $\rho_i=\bb{R}_{\geq 0} \?{e}_i$.  Assume the elements of $F$ are labelled cyclically by $\{1,\ldots,c\}$ for some $c\geq 3$, and let $\sigma_{i,i+1}$ denote the cone bounded by $\rho_i$ and $\rho_{i+1}$.  By standard toric geometry, $\pi_2(e_{i-1})+\pi_2(e_{i+1})+[\?{D}_i]^2\pi_2(e_i)=0$.  So by the definition of a bending parameter, 
\begin{align*}
    \varphi(\?{e}_{i+1})&=-\varphi(\?{e}_{i-1})-[\?{D}_i]^2\varphi(\?{e}_i)+\f{b}_iu_i(\?{e}_{i+1})
\end{align*}
where $u_i$ is the primitive element of $\f{d}_i^{\perp}\subset \?{M}$ which is positive on $\?{e}_{i+1}$.  So
\begin{align*}
    \f{b}_i=e_{i-1}+e_{i+1}+[\?{D}_i]^2 e_i = \kappa^{-1} \circ \eta_{\SSS}^* ([\?{D}_i]),
\end{align*}
as desired.

Now, letting $\?{e}_i:=\pi_2(e_i)$ for $i\in I_{\uf}$ as well, the scattering diagram $\?{\f{D}}_0$ is defined by 
\begin{align*}
    \?{\f{D}}_0:=\{(\bb{R} \?{e}_i, 1+z^{(\?{e}_i,\varphi'(\?{e}_i)-[E_i]}))|i \in I_{\uf}\}.
\end{align*}
We wish to show that this agrees with $\f{D}_{\In}^{\SSS}$ as in \eqref{Din}.  When expressing Looijenga pairs in terms of cluster varieties we choose each $\pi_1(e_i)$ to be primitive, and in rank $2$ it is clear that $\bb{R} e_i=\pi_1(e_i)^{\perp}$, so we need only check that $e_i=\varphi'(\?{e}_i)-\kappa^{-1}([E_i])$.  By construction, $\?{e}_i=\?{e}_j$ for some $j\in F$, and $\varphi'(\?{e}_j)=e_j$, so the claim is that $\kappa(e_j-e_i)=[E_i]$, which is immediate from the construction of $\kappa$.  Thus, $\f{D}^{\SSS}=\?{\f{D}}$ and the result follows.
\end{proof}

It follows that our main results, Theorems \ref{MainNaive} and \ref{MainThm}, will apply in particular to the theta functions of \cite{GHK1}.

\section{Theta functions, tropical curves, and log GW invariants of toric varieties}\label{Review}
\subsection{Tropical description of theta functions}\label{ScatTheta}

In this subsection we summarize some results from \cite{Man3} which express the theta functions in terms of certain counts of tropical curves.  We work with a fixed seed $\SSS$.

\subsubsection{The module of log derivations}\label{logDerivations}

Our setup for describing scattering diagrams in \S \ref{ScatterIntro} is a simplified version of the more general setup used in \cite{Man3}.  We briefly explain how to translate between the two setups so that we can apply the results of \cite{Man3} to the setup here.  See also \cite[Examples 2.1(i), 2.7(i), and 3.4(i)]{Man3} for specializations of the general setup there to our situation, as well as \cite[Rmk. 2.3(i)]{Man3} for an explanation of an additional difference between the setups (roughly, supports of walls here are equal to $\pi_2$ of the supports of walls there).

Recall that $\?{M}:=\pi_1(N)\subset M$ and $R:=\kk[K_2\cap P]$ where $K_2=\ker(\pi_2)$.  As in \cite[\S 1.1]{GPS}, consider the module of log derivations $\Theta_K(N^{\oplus})$ defined by
\begin{align*}
    \Theta_K(N^{\oplus}):=\kk[N^{\oplus}]\otimes_{\bb{Z}}\?{M}
\end{align*}
with action on $\kk[N^{\oplus}]$ via $R$-derivations defined by
\begin{align*}
    f\otimes m(z^n):=f\langle n,m\rangle z^n.
\end{align*}
We will write $f\otimes m$ as $f\partial_m$.  $\Theta_K(N^{\oplus})$ forms a Lie algebra with bracket $[a,b]:=ab-ba$, where multiplication means composition of derivations.  In particular,
\begin{align*}
    [z^{n_1}\partial_{m_1},z^{n_2}\partial_{m_2}]=z^{n_1+n_2}\partial_{\langle n_2,m_1\rangle m_2-\langle n_1,m_2\rangle m_1}.
\end{align*}
For each $n\in N^+$, let $\f{h}_n$ be the submodule of $\Theta_K(N^{\oplus})$ spanned by the element $z^n\partial_{\pi_1(n)}$.  One easily checks that $\f{h}:=\bigoplus_{n\in N^+} \f{h}_n$ is a Lie subalgebra of $\Theta_K(N^{\oplus})$.  Let $\wh{\f{h}}$ denote the completion of $\f{h}$ associated to the $N^+$-grading.  For $n\in N^{+}$ primitive, let $\f{h}_n^{\parallel}$ denote the Lie subalgebra of $\wh{\f{h}}$ spanned topologically by elements of $\f{h}_{kn}$ for $k\in \bb{Z}_{>0}$.

Now, in the setup of \cite{Man3}, our walls would be written as $(m,\f{d},g)$, where $m$ is an element of $\?{M}$ up to positive scaling such that the support $\f{d}$ is contained in $m^{\perp}\subset \?{N}_{\bb{R}}$, and $g$ is an element of $\f{h}_n^{\parallel}$ for some primitive $n\in m^{\perp}\subset N$ (recall that $m\in \?{M}\subset M$, so we can view $m^{\perp}$ as living in $\?{N}_{\bb{R}}$ or $N$).  More precisely, consider a wall given in our setup by $(\f{d},f\in \kk\llb n\rrb)$, and let $u$ be a primitive element of $\f{d}^{\perp}\cap \?{M}$.  Then in the setup of \cite{Man3}, this wall would be expressed as $(u,\f{d},\log(f)\partial_{u})$, cf. \cite[Ex. 2.8(i)]{Man3}.  The wall-crossing automorphism $\theta_{\f{d}}$ as in \eqref{WallCross} is then realized as the action of $\exp \left[\sign \langle u,-\gamma'(t)\rangle \ad_{\log(f)\partial_u}\right]$.

In particular, recall the scattering diagram $\f{D}_{\In}:=\left\{\left(\pi_1(e_i)^{\perp},(1+z^{e_i})^{|\pi_1(e_i)|}\right) : i\in I_{\uf}\right\}$ from \eqref{Din}.  In the setup of \cite{Man3}, this could be written as
\begin{align}\label{DinMan3}
    \f{D}_{\In}=\left\{\left( \pi_1(e_i),\pi_1(e_i)^{\perp},g_i:=\sum_{w=1}^{\infty} w \frac{(-1)^{w+1}}{w^2} z^{we_i} \partial_{\pi_1(e_i)}\right) : i \in I_{\uf} \right\}.
\end{align}
Note that $\pi_1(e_i)=d_i\omega(e_i,\cdot)$, so $g_i$ is in fact in $\wh{\f{h}}$.  Denoting $a_{iw}:=w\frac{(-1)^{w+1}}{w^2}$ (pulling a $w$-factor in front like this will be convenient in \S \ref{DegenS}), we can rewrite $g_i$ as $g_i=\sum_{w=1}^{\infty} g_{iw}$ where
\begin{align}\label{gi}
    g_{iw}:= a_{iw} z^{we_i} \partial_{\pi_1(e_i)}.
\end{align}

\subsubsection{Tropical Gromov-Witten invariants}\label{TropGW}
We now recall some background on tropical curves.  It is convenient to work in greater generality than we will need.  We begin by recalling the definition of a (genus $0$) tropical curve, following the conventions of \cite[\S 3.1]{Man3} or \cite[\S 2]{MRud} which built on the definitions of tropical curves in \cite[Def. 2.2]{Mi} and \cite[\S 1]{NS}.

 Let $\?{\Gamma}$ be the topological realization of a connected finite tree without bivalent vertices, and denote the complement of the $1$-valent vertices by $\Gamma$.  Let  $\Gamma^{[0]}$, $\Gamma^{[1]}$, and $\Gamma^{[1]}_{\infty}$ denote the vertices, edges, and unbounded edges of $\Gamma$, respectively.   We equip $\Gamma$ with a weight function $\ww:\Gamma^{[1]}\rar \bb{Z}_{\geq 0}$ such that if $\ww(E)=0$, then $E\in \Gamma^{[1]}_{\infty}$.  A marking of $\Gamma$ is a bijection $\epsilon:J\rar \Gamma^{[1]}_{\infty}$ for some index set $J$.  We denote $E_j:=\epsilon(j)$.  Let $J^{\circ}\subset J$ denote the set of $j\in J$ for which $\ww(E_j)=0$, and let $J':=J\setminus J^{\circ}$.

\begin{dfn}\label{TropCurveDfn}
A (genus $0$) parameterized tropical curve $(\Gamma,h)$ (in $\?{N}_{\bb{R}}$) is data $\Gamma$, $\ww$, $\epsilon$ as above (the weight and marking are suppressed in the notation), along with a continuous map $h:\Gamma\rar \?{N}_{\bb{R}}$ such that
\begin{itemize}
\item For each $E\in \Gamma^{[1]}$ with $\ww(E)>0$, $h|_E$ is a proper embedding into an affine line with rational slope.  If $\ww(E)=0$, then $h(E)$ is a point.
\item \textbf{The balancing condition}: For any edge $E\ni V$ with $\ww(E)>0$, denote by $u_{(V,E)}$ the primitive integral vector emanating from $h(V)$ into $h(E)$.  Then
\begin{align}\label{balance}
\sum_{E\ni V} \ww(E) u_{(V,E)} =0.
\end{align}
\end{itemize}
For non-compact edges $E_i\ni V$, we may denote $u_{(V,E_i)}$ as simply $u_{E_i}$.  Similarly, we may write $u_E$ in place of $u_{(V,E)}$ if $V$ is clear or if the direction only matters up to scaling.

 An isomorphism of marked parameterized tropical curves $(\Gamma,h)$ and $(\Gamma',h')$ is a homeomorphism $\Phi:\Gamma\rar \Gamma'$ respecting the weights and markings such that $h=h'\circ \Phi$.  A {\bf tropical curve} is then defined to be an isomorphism class of parameterized marked tropical curves.  We will let $(\Gamma,h)$ denote the isomorphism class it represents, and we will often abbreviate this as simply $\Gamma$ or $h$.
 
In \cite[\S 3.1]{Man3}, a \textbf{tropical disk} is defined in the same way except that there is an edge $E_{\out}\in \Gamma^{[1]}_{\infty}$ of possibly positive weight which is contracted by $h$.  The associated label in $J$ is viewed as being in $J^{\circ}$, not $J'$.  Let $V_{\out}$ denote the vertex of $E_{\out}$.  The balancing condition is still required to hold at the vertex $V_{\out}$ of $E_{\out}$ for some uniquely determined primitive vector $u_{(V_{\out},E_{\out})}\in \?{N}$ (or for $u_{(V_{\out},E_{\out})}=0$ if $\ww(E_{\out})=0$, in which case the tropical disk can be viewed as a tropical curve).
\end{dfn}

Let $\val(V)$ denote the number of edges containing $V$ (the valence of $V$).  
Let $\Flags(\Gamma)$ denote the set of flags $(V,E)$, $V\in E$, of $\Gamma$.  The {\bf type} of a tropical curve or disk is the data of $\Gamma$, $\ww$, and $\epsilon$, along with the data of the map $u:\Flags(\Gamma)\rar \?{N}$, $(V,E)\mapsto u_{(V,E)}$.  The \textbf{degree} $\Delta$ of a tropical curve or disk is the data of $J$ along with the corresponding map $\Delta:J\rar \?{N}$, $\Delta(j)=\ww(E_j)u_{E_j}$.  In the case of a tropical disk, we say that $\Delta$ also remembers which edge is the special one $E_{\out}$ which is contracted despite possibly having positive weight.

\begin{dfn}\label{Constraints}
An affine constraint $\AAA$ is a tuple $(\A_j)_{j\in J}$ of affine subspaces of $\?{N}_{\bb{R}}$ with rational slope, each equipped with a weight\footnote{Most of \cite{MRud} treats these weights as being $1$, but a generalization allowing for other weights appears in \cite[Def. 4.17]{MRud}.  Multiplying $\ww(\A_j)$ by $k$ just has the effect of multiplying the multiplicities of our tropical curves by $k$.} $\ww(\A_j)\in \bb{Z}_{>0}$.  A tropical curve or disk $(\Gamma,h)$ matches the constraint $\AAA$ if $h(E_j)\subset \A_j$ for all $j\in J$.  
Now consider a map $\Psi:J^{\circ}\rar \bb{Z}_{\geq 0}$, denoting $s_j:=\Psi(j)$.  We say $(\Gamma,h)$ satisfies the $\psi$-class conditions $\Psi$ if \begin{align}\label{Psidfn}
    \val(V)-3\geq \sum_{\substack{j\in J^{\circ}\\ E_j\ni V}} s_j
\end{align}
for each vertex $V\in \Gamma^{[0]}$.  Let $\langle V\rangle$ denote the multinomial coefficient
\begin{align}\label{V}
    \langle V \rangle:= \frac{(\val(V)-3)!}{\prod_{\substack{j\in J^{\circ}\\ E_j\ni V}} s_j!}
\end{align}
and denote $\langle \Gamma \rangle:=\prod_{V\in \Gamma^{[0]}} \langle V\rangle$.
\end{dfn}

Let $\f{T}_{0,\Delta}(\AAA,\Psi)$ denote the space of genus $0$ degree $\Delta$ tropical curves or disks which match the constraint $\AAA$ and satisfy a $\psi$-class condition $\Psi$.  By \cite[Lem 2.14]{MRud} (which is stated for tropical curves but easily extends to allow for tropical disks), this space is finite (and \eqref{Psidfn} is an equality for each $V$) for generic translates of the $\A_j$'s whenever 
\begin{align}\label{RigidDim}
 \sum_{j\in J} \codim(\A_j) + \sum_{j\in J^{\circ}} s_j= |J|+r-3.    
\end{align}
Under these conditions (i.e., when \eqref{RigidDim} is satisfied and the $\A_i$'s are chosen generically among their translations), $\f{T}_{0,\Delta}(\AAA,\Psi)$ and the tropical curves/disks it contains are called \textbf{rigid}.  So assuming rigidity of $\f{T}_{0,\Delta}(\AAA,\Psi)$, after describing how to assign a ``multiplicity''  $\Mult(\Gamma)\in \bb{Z}_{\geq 1}$ to each $\Gamma\in\f{T}_{0,\Delta}(\AAA,\Psi)$,
 one can define tropical Gromov-Witten invariants as follows:\footnote{Counts of tropical curves weighted by their multiplicities were first considered (and related to holomorphic curve counts) in \cite{Mi} and \cite{NS}.  Contributions of $\psi$-class conditions in certain cases were considered in \cite{MR,GrP2}.  The setup here follows \cite{MRud}.}
\begin{dfn}
\begin{align}\label{GWtropDfn}
\GW_{0,\Delta}^{\trop}(\AAA,\Psi):=\sum_{(\Gamma,h)\in \f{T}_{0,\Delta}(\AAA,\Psi)} \langle \Gamma \rangle \Mult(\Gamma).
\end{align}
\end{dfn}

For the cases we care about, $s_j$ will be nonzero for only one $j$, and so since rigidity ensures that \eqref{Psidfn} is always an equality, $\langle \Gamma \rangle$ will always be $1$.  We note that \eqref{GWtropDfn} is independent of the generic translates of the $\A_j$'s for tropical curves, but for tropical disks we will have to specify $\A_{\out}$ more precisely.

\subsubsection{Multiplicities of tropical curves and disks}\label{MultGWtrop}

The definition of $\Mult(\Gamma)$ used in \cite[Lem/Def 2.16]{MRud} (due to \cite{NS} when there are no $\psi$-classes) is given in a form which is impractical for the applications we consider here.  This motivated the paper \cite{MRudMult}, which shows that the same multiplicities can alternatively be computed as follows:

Consider a rigid $(\Gamma,h)\in \f{T}_{0,\Delta}(\AAA,\Psi)$.  For each $j\in J$, let $\alpha_j$ denote an index-$\ww(\A_j)$ element (unique up to sign) of $\Lambda^{\codim \A_j} \?{M} \subset \Lambda^* \?{M}$ whose kernel is parallel to $\A_j$, i.e., the contraction\footnote{Given a lattice $L$ with dual $L^*$, recall that an element $\alpha\in \Lambda^k L^*$ can be viewed as a multilinear function $\alpha:L^k\rar \bb{Z}$.  Then for $n\in L$, the contraction $\iota_n(\alpha)\in \Lambda^{k-1}L^*$ is the element corresponding to the multilinear function $\iota_n(\alpha):L^{k-1}\rar \bb{Z}$, $(n_1,\ldots,n_{k-1})\mapsto \alpha(n,n_1,\ldots,n_{k-1})$.} $\iota_n(\alpha_j)=0$ if and only if $n$ is parallel to $\A_j$.  We pick a flow on $\Gamma$ by choosing one vertex to serve as the sink.   Using this flow, we will recursively associate an element $\omega_E$ of $\bb{Z}[\?{N}]\otimes \Lambda^* \?{M}$ (determined up to sign) to every edge $E$ of $\Gamma$.  For each $j\in J$, we associate the element $$\omega_{E_j}:=z^{\ww(E_j)u_{E_j}}\otimes \alpha_j$$ to the edge $E_j$.  Now consider a vertex $V\neq V_0$ with $E_1,\ldots,E_k$ flowing into $V$ and $E_V$ the unique edge flowing out of $V$.  If $\omega_{E_i} = z^{n_i}\otimes \alpha_i$ is the element associated to $E_i$, $i=1,\ldots,k$, we define the element associated to $E_V$ to be 
\begin{align}\label{lk}
    \omega_{E_V} := {\bf L}_k(z^{n_1}\otimes \alpha_1,\ldots,z^{n_k} \otimes \alpha_k) := z^{n_V} \otimes \iota_{n_V}(\alpha_1 \wedge \cdots \wedge \alpha_k),
\end{align}
where $n_V:=n_1+\ldots+n_k$.

Finally, if $E_1,\ldots,E_s$ are the edges containing $V_0$, define \begin{align*}
    \omega_0:=\omega_{E_1} \cdots  \omega_{E_s} \in 1\otimes \Lambda^{r} \?{M} \subset \bb{Z}[\?{N}]\otimes \Lambda^* \?{M}.
\end{align*}
Here, balancing ensures that the first factor is $1$ (the exponents of the monomial terms cancel out), and rigidity ensures that the wedge product is of top degree.  Then \cite[Thm. 1.2]{MRudMult} says that the multiplicity of $\Gamma$ considered in \cite{MRud} can be computed as
\begin{align}\label{MultDfn}
    \Mult(\Gamma)=|\omega_0|
\end{align}
by which we mean the index of $\omega_0$ in $\Lambda^r \?{M}$.

Finally, while \cite{MRudMult} focused on tropical curves, we note that the above recipe yields a well-defined number $\Mult(\Gamma)$ when $\Gamma$ is a rigid tropical disk as well.

\subsubsection{Theta functions from tropical curve counts}\label{ThetaFromTrop}

We next describe the specific tropical degrees and conditions used in the main theorems of \cite{Man3}, applied to the scattering diagram $\f{D}:=\scat(\f{D}_{\In})$ for $\f{D}_{\In}$ the initial scattering diagram \eqref{DinMan3} associated to a seed $\SSS=(N,I,E:=\{e_i\}_{i\in I},F,B)$.

Let $\ww:=(\ww_i)_{i\in I_{\uf}}$ be a tuple of weight vectors $\ww_i:=(w_{i1},\ldots,w_{il_i})$ with $w_{i1} \leq \ldots \leq w_{il_i}$, $w_{ij}\in \bb{Z}_{>0}$.  For $\Sigma_{l_i}$ denoting the group of permutations of $\{1,\ldots,l_i\}$, let \[\Aut(\ww)\subset \prod_{i\in I_{\uf}} \Sigma_{l_i}\] be the group of automorphisms of the second indices of the weights $\ww_i$ which act trivially on $\ww$.

Let $\pp$ be an $s$-tuple $(p_1,\ldots,p_s)$ of elements of $\?{N}$.  For $n\in N$ and $\varphi:\?{N}_{\bb{R}}\rar N_{\bb{R}}$ as in \eqref{varphi}, let $\s{W}_{\pp}(n)$ denote the set of weight vectors $\ww$ such that
\begin{align*}
    \sum_{i\in I_{\uf}} \sum_{j=1}^{l_i} w_{ij}e_i +\sum_{k=1}^s \varphi(p_k)= n.
\end{align*}
Let \begin{align*}
    J_{\ww,\pp}:=\{(i,j)|i\in I_{\uf},j=1,\ldots,l_i\} \cup \{1,\ldots,s\}\cup \{\out,\infty\}.
\end{align*}
We will often write the pairs $(i,j)$ as simply $ij$.

For $\ww\in \s{W}_{\pp}(n)$, let $\Delta_{\ww,\pp}$ denote the degree $$\Delta_{\ww,\pp}:J_{\ww,\pp}\rar \?{N}$$ with $\Delta_{\ww,\pp}((i,j)):=w_{ij}\pi_2(e_i)$, $\Delta_{\ww,\pp}(k)=p_k$ for $k=1,\ldots,s$, $\Delta_{\ww,\pp}(\out)=-\pi_2(n)$, and $\Delta_{\ww,\pp}(\infty)=0$.  We view $\Delta_{\ww,\pp}$ as the degree of a tropical disk, with $E_{\out}$ being the special contracted edge.  We have $$J^{\circ}_{\ww,\pp}=\Delta_{\ww,\pp}^{-1}(0)\cup \{\out\}=\{1,\ldots,s,\out,\infty\}.$$

Given a generic point $Q\in \?{N}_{\bb{R}}$ and $\ww\in \s{W}_{\pp}(n)$, we define the incidence conditions $\AAA_{\ww,\pp,Q}$ as follows:
Each $\A_{ij}$ is taken to be a generic translate of $\f{d}_i=\pi_1(e_i)^{\perp}\subset \?{N}_{\bb{R}}$ with $\ww(\A_{ij})=|\pi_1(e_i)|$.  We take $\A_k:=\?{N}_{\bb{R}}$ for each $k=1,\ldots,s$ and for $k=\infty$ (i.e., the incidence conditions on the $E_k$'s are trivial), and we take $\A_{\out}=Q$.  For our $\psi$-class conditions, we define $\Psi_{\ww,\pp}:J_{\ww,\pp}^{\circ}\rar \bb{Z}_{\geq 0}$ by $\Psi_{\ww,\pp}(\out)=s-1$ and $\Psi_{\ww,\pp}(k)=0$ for every other $k\in J_{\ww,\pp}^{\circ}$.

\begin{rmk}
In the setup of \cite{Man3}, the marked point $\infty$ is not included, but $\Psi_{\ww,\pp}(\out)$ is taken to be $s-2$ instead of $s-1$.  This change in the $\psi$-class condition forces the valence of the vertex $V_0\in E_{\out}$ to be higher by $1$ here, thus forcing the extra contracted edge $E_{\infty}$ to contain $V_0$, so there is an obvious bijection between the tropical curves/disks in the two setups.  This modification of the tropical data corresponds to the geometric modification discussed in Remark \ref{logFCA}, and it allows us to avoid treating the $s=1$ case separately.
\end{rmk}

With these conditions and $\ww \in \s{W}_{\pp}(n)$, $\f{T}_{0,\Delta_{\ww,\pp}}(\AAA_{\ww,\pp,Q},\Psi_{\ww,\pp})$ is finite, so we can count its elements using the multiplicities $\Mult$ introduced above, yielding numbers $\GW^{\trop}_{0,\Delta_{\ww,\pp}}(\AAA_{\ww,\pp,Q},\Psi_{\ww,\pp})$.  Alternatively, \cite[\S 3.1.2]{Man3} defines slightly different multiplicities $\mult$ (explained below), and with these one defines 
\begin{align*}
    \N^{\trop}_{\ww,\pp}(Q):=\sum_{\Gamma \in \f{T}_{0,\Delta_{\ww,\pp}}(\AAA_{\ww,\pp,Q},\Psi_{\ww,\pp})} \mult(\Gamma) \in \wh{A}.
\end{align*}

For $\Gamma \in \f{T}_{0,\Delta_{\ww,\pp}}(\AAA_{\ww,\pp,Q},\Psi_{\ww,\pp})$, we define the multiplicity $\mult(\Gamma)$ by modifying the definition of $\Mult(\Gamma)$ as follows: For the sink $V_0$ of our flow, we use the vertex $V_{\out}$ contained in $E_{\out}$.  Then in the recursive construction, we associate to each edge $E$ an element $\omega_E$ of the Lie algebra $A\oplus \f{h}$ for $A$ as in \S \ref{ScatterIntro} and $\f{h}$ as in  \S \ref{logDerivations}.  The bracket is given by 
\begin{align}\label{AhBracket}
    [a_1 + h_1,a_2 +h_2]=[a_1,a_2]+h_1\cdot a_2 - h_2\cdot a_2 + [h_1,h_2],
\end{align} where $\cdot$ is the derivation action of $\f{h}$ on $A$.  

To $E_{\out}$ and $E_{\infty}$ we associate the element $1\in A$.  To $E_{ij}$ we associate the element $g_{iw_{ij}}$ as defined in \eqref{gi}.  To $E_k$ for $k=1,\ldots,s$ we associate $z^{\varphi(p_k)}\in A$.  Then at a vertex $V\neq V_0$ (necessarily trivalent since the affine constraints are translated generically, cf. \cite[Lem. 2.14]{MRud}), instead of applying ${\bf L}_2$ as in \eqref{lk}, we simply take the Lie bracket of the two incoming elements (for a certain choice of ordering) in order to produce the outgoing $\omega_{E_V}$.  Finally, for each $E\ni V_0$, $\omega_E$ is in fact an element of $A$, and the product $\prod_{E\ni V} \omega_E$ in $A$ is, up to sign, equal to $\mult(\Gamma)$.  This element $\mult(\Gamma)\in A$ has the form $ka_{\ww}z^{n_{\out}}$ for some nonzero integer $k$,
\begin{align*}
    a_{\ww}:=\prod_{ij} a_{iw_{ij}},
\end{align*} and $n_{\out}:=\sum_{ij} w_{ij}e_i + \sum_{k=1}^s \varphi(p_k)$, i.e., $n_{\out}$ is the element of $N$ such that $\ww\in \s{W}_{\pp}(n_{\out})$.   As explained in \cite[Ex. 3.4(i)]{Man3}, the correct sign of $\mult(\Gamma)$ is the one for which $k$ is positive.

Now, in the definition of $\Mult(\Gamma)$, the elements $g_{iw_{ij}}=a_{iw_{ij}}z^{w_{ij}e_i}\partial_{\pi_1(e_i)}$ above would have instead been $z^{w_{ij}\pi_2(e_i)}\partial_{\pi_1(e_i)}$, while the elements associated to $E_k$ for $k=1,\ldots,s$ would have been $z^{p_k}$ instead of $z^{\varphi(p_k)}$.  Then in the computation of $\Mult$, the element $\omega_{E_{\out}}$ is $z^{-\pi_2(n_{\out})}$ times a primitive element of $\Lambda^{\topp} M$.  The element $\omega_{E_{\infty}}$ is still $1$.  Note that  $\pi_2:N\rar \?{N}$ induces a map $\pi_2:A\oplus \f{h}\rar \bb{Z}[\?{N}]\otimes \Lambda^* \?{M}$ (i.e., applying $\pi_2$ to the exponents), and for ${\bf L}_2$ as in \eqref{lk} and any $a,b\in A\oplus \f{h}$, we have ${\bf L}_2(\pi_2(a),\pi_2(b))=\pi_2([a,b])$ for $[\cdot,\cdot]$ as in \eqref{AhBracket}.  One now checks that
\begin{align}\label{mM1}
    \mult(\Gamma)=a_{\ww}\Mult(\Gamma)z^{n_{\out}}.
\end{align}

\subsubsection{Theta functions and scattering diagrams in terms of tropical invariants}

We are now ready to state \cite{Man3}'s result expressing theta functions in terms of tropical disk counts.

\begin{lem}[\cite{Man3} Thm. 3.9 and Prop 2.15]\label{TropFrob}
For $\pp:=(p_i)_{i=1,\ldots,s}$ an $s$-tuple of elements of $\?{N}$, let $\alpha(\pp;p)$ denote the $\vartheta_p$-coefficient of $\prod_{i=1}^s \vartheta_{p_i}$.  Then
\begin{align}\label{TropFrobEqnQuot}
    \alpha(\pp;p)z^{\varphi(p)}=\sum_{n\in K_2^+} \sum_{\ww\in \s{W}_{\pp}(\varphi(p)+n)} \frac{\N_{\ww,\pp}^{\trop}(Q)}{|\Aut(\ww)|}
\end{align}
for $Q$ chosen to be sufficiently close to the ray through $p$ (or for $Q$ anywhere if $p=0$).
\end{lem}

\begin{rmk}\label{Dints}
Using \eqref{mM1}, note that we can rewrite \eqref{TropFrobEqnQuot} as
\begin{align}\label{TropFrobEqnQuot2}
        \alpha(\pp;p)=\sum_{n\in K_2^+} \sum_{\ww\in \s{W}_{\pp}(\varphi(p)+n)} a_{\ww}\frac{\GW^{\trop}_{0,\Delta_{\ww,\pp}}(\AAA_{\ww,\pp,Q},\Psi_{\ww,\pp})}{|\Aut(\ww)|}z^{n}.
\end{align}
Given data $\pp$, $p$, $n$, and $\ww\in \s{W}_{\pp}(\varphi(p)+n)$ as in \eqref{TropFrobEqnQuot2}, it follows immediately from the definition of $\Delta_{\ww,\pp}$ and the description of $\kappa$ in Theorem \ref{kappa} that $\kappa(n)\in N_1(Y_S,\bb{Z})$ is the class of a curve having intersection number $\sum_{j=1}^{l_i} w_{ij}$ with $[E_i]$ for each $i\in I_{\uf}$, plus for each ray $\rho$ of $\Sigma$ having intersecting number with $D_{\rho}$ equal to $\sum_{p_i\in \rho} |p_i|,$ plus $-|p|$ if $p\in \rho$.
\end{rmk}

\begin{cor}\label{AmplePolynomial}
If $(Y_{\SSS,\Sigma}, D_{\SSS,\Sigma})$ supports an ample divisor on its boundary, then the theta functions in fact generate an algebra over $\kk[\NE(Y_{\SSS,\Sigma})]$ (as opposed to over some formal completion of this).
\end{cor}
\begin{proof}
By possibly replacing $\Sigma$ by a refinement, we have that $D_{\SSS,\Sigma}$ supports an effective ample divisor $H$.  After possibly replacing $H$ with $D_{\SSS,\Sigma}+kH$ for $k$ a sufficiently large integer, we can always assume that $H$ is of the form $\sum_{\rho\in \Sigma^{[1]}} a_{\rho}D_{\rho}$ with each $a_{\rho}$ in $\bb{Z}_{\geq 1}$.  By Lemma \ref{NECoeff}, the coefficient of $\vartheta_p$ in the theta function expansion of the product $\prod_{i=1}^s \vartheta_{p_i}$ is an (a priori formal) sum of monomials whose exponents correspond to elements of the Mori cone $\NE(Y_{\SSS,\Sigma})$, hence have positive intersection with $H$.  Furthermore, Remark \ref{Dints} ensures that the intersections of these classes with $H$ are determined by $\pp$ and $p$, and for fixed $\pp$ we see that there are only finitely many possibilities for $p$ for which this intersection number with $H$ will in fact be positive.  Thus, there are only finitely many $p$'s for which the $\vartheta_p$-coefficient can be nonzero.  

So now it suffices to show that for fixed $\pp$ and $p$, there are only finitely many curve classes with the associated intersection number $d$ with $H$.  Note that there always exists a rational polyhedral cone $\Xi\subset N_1(Y_{\SSS,\Sigma})$ such that $\NE(Y_{\SSS,\Sigma})\subset \Xi$ and such that $H$ is positive on $\Xi\setminus \{0\}$.\footnote{In fact, as the author has learned from Sean Keel,  $D_{\SSS,\Sigma}$ supporting an ample divisor implies that $(Y_{\SSS,\Sigma}, D_{\SSS,\Sigma})$ is log Fano, hence that $\NE(Y_{\SSS,\Sigma})$ itself is rational polyhedral by \cite[Cor. 1.3.2]{BCHM}.}  Since $\Xi$ is finitely generated, it clearly contains only finitely many elements whose pairing with $H$ is $d$, and so the same holds for $\NE(Y_{\SSS,\Sigma})\subset \Xi$, as desired.
\end{proof}

\subsection{Log Gromov-Witten invariants}\label{SectionLogGW}

Let us recall and extend the setup from \S \ref{FSC}.  Given a smooth log pair $Y^{\dagger}=(Y,D)$, we can consider the integral points $B(\bb{Z})$ of the cone over the dual intersection complex of some $\wt{D}$.  In particular, when $(Y,D)$ is  $(Y_{\SSS,\Sigma},D_{\SSS,\Sigma})$ or $(\TV(\wt{\Sigma}),\partial \TV(\wt{\Sigma}))$ for $\wt{\Sigma}$ a complete fan in $\?{N}$, the set $B(\bb{Z})$ is identified with $\?{N}$.

Recall that a tropical degree is a map $\Delta:J\rar B(\bb{Z})$ for some finite index-set $J$.  Let $J^{\circ}:=\Delta^{-1}(0)$.  Let $\Sigma$ be the fan in $B$ for $Y^{\dagger}$.  Let $\eta:(\wt{Y},\wt{D})\rar (Y,D)$ denote a toric blowup corresponding to a refinement $\wt{\Sigma}$ of $\Sigma$.  We can assume that $\wt{Y}$ is projective and that $\Delta$ is a {\bf torically transverse degree} for $\wt{\Sigma}$, meaning that $\Delta(j)$ is contained in a ray $\rho\in \wt{\Sigma}$ for each $j\in J'$.  Recall that $D_{\Delta(j)}\subset \wt{D}$ denotes the corresponding boundary component. 

For $\beta\in \NE(\wt{Y})$, we consider the moduli stack $\s{M}_{0,\Delta}^{\log}(\wt{Y}^{\dagger},\beta)$ as in \S \ref{FSC}.   This moduli stack was constructed and shown to be an algebraic stack in \cite{GSlog} and \cite{AC}.  Furthermore, \cite[Thm. 0.3]{GSlog} says that $\s{M}_{0,\Delta}^{\log}(\wt{Y}^{\dagger},\beta)$ is equipped with a virtual fundamental class $[\s{M}^{\log}_{0,\Delta_{\pp}}(\wt{Y}^{\dagger},\beta)]^{\vir}$ of virtual dimension $$\vdim(\s{M}^{\log}_{0,\Delta_{\pp}}(\wt{Y}^{\dagger},\beta))=\dim(\wt{Y})-(K_Y+D)\cdot \beta+|J|-3$$ which satisfies the standard properties of virtual fundamental classes, thus allowing for the construction of log Gromov-Witten invariants:
\begin{dfn}
Consider the evaluation maps $\ev_i$ for each $i\in J$ and $\psi$-classes $\psi_j$ as in \eqref{psic1}.  Define a map $Z:J\rar A^*(\wt{Y})$, denoting $Z_j:=Z(j)$.  Assume each $Z_j$ is the class of a regularly embedded closed subvariety of $Y$.  Let $\Psi$ be a map from $J^{\circ}$ to $\bb{Z}_{\geq 0}$, and denote $s_j:=\Psi(j)$.  Then one defines\footnote{More precisely, let $\gamma :=  \left(\prod_{j}  \psi_j^{s_j} \right) \cap \left(\prod_{i} [Z_i] \right) \cap [\s{M}^{\log}_{0,\Delta_{\pp}}(\wt{Y}^{\dagger},\beta)]^{\vir} \in A^*(\s{M}^{\log}_{0,\Delta_{\pp}}(\wt{Y}^{\dagger},\beta))$, where the cap products are interpreted using generalized Gysin maps for the cycles $[Z_i]$ and capping with Chern classes for the classes $\psi_j^{s_j}$.  Then $\GW^{\log}_{0,Y^{\dagger},\Delta, \beta}(Z,\Psi) := \deg(\gamma)\in \bb{Q}$.  Here, $\deg(\gamma):=0$ unless  $\vdim(\s{M}^{\log}_{0,\Delta_{\pp}}(\wt{Y}^{\dagger},\beta))-\sum_i \codim(Z_i) - \sum_j s_j=0$.  Cf. \cite[Def. 3.7 and \S A]{MRud} for more on the technical details of this definition.} \textbf{log Gromov-Witten invariants} by
\begin{align}\label{GWdef}
\GW^{\log}_{0,Y^{\dagger},\Delta, \beta}(Z,\Psi) := \int_{[\s{M}^{\log}_{0,\Delta}(\wt{Y}^{\dagger},\beta)]^{\vir}} \left(\bigcup_{j\in J^{\circ}} \psi_j^{s_j}\right) \cup \left(\bigcup_{j\in J} \ev_j^*(Z_j)\right).
\end{align}
\end{dfn}
The main result of \cite{AW} ensures that the log Gromov-Witten invariants are independent of the choice of toric blowup $\wt{Y}$ of $Y$.

Suppose $(\wt{Y},\wt{D})$ above is a nonsingular complete toric variety with cocharacter lattice $\?{N}$ and fan $\wt{\Sigma}$, so $B(\bb{Z})=\?{N}$.  Consider an affine linear subspace $\A\subseteq \?{N}_{\bb{R}}$ with rational slope.  Given a point $x$ in the big torus orbit of $Y$, we obtain a subvariety $Z_{\A,x}$ as follows: Let $\A^{\perp}$ denote the $m\in \?{M}$ which pair to $0$ with the tangent directions to $\A$.  Let $z^m(x)$ denote $z^m$ evaluated at $x$.  Then $Z_{\A,x}$ is the subvariety corresponding to the ideal sheaf $\langle z^m-z^m(x)|m\in \A^{\perp}\rangle$.  In particular, when $\A\subset\?{N}_{\bb{R}}$ is just a point, $Z_{\A,x}=x$.  Note that $\dim_{\kk}(Z_{\A,x})=\dim_{\bb{R}}(\A)$.  We denote the corresponding Chow class (independent of $x$) by $[Z_{\A}]$.

Now suppose we have data $\?{N},\Delta,\AAA,\Psi$ as in \S \ref{TropGW}.  Define $Z_{\AAA}:J\rar A_*(\wt{Y})$ by $Z_{\AAA}(j)=[Z_{\A_j}]$.  Let $[\Delta]\in \NE(\wt{Y},\bb{Z})$ be the unique curve class such that $$[\Delta].[D_{\rho}]=\sum_{\substack{j\in J \\ \Delta(j)\in \rho}} |\Delta(j)|$$ for each ray $\rho$ of $\wt{\Sigma}$.  We define
\begin{align}\label{GWtoric}
\GW_{0,\Delta}^{\log}(\AAA,\Psi):=\GW^{\log}_{0,\TV(\wt{\Sigma})^{\dagger},\Delta,[\Delta]}(Z_{\AAA},\Psi).
\end{align}

The genus $0$ case of \cite[Thm 4.15]{MRud}, or alternatively \cite[Cor 5.2]{AGr}, states the following ($\Delta$ here is assumed to be a tropical curve degree, not a tropical disk degree):
\begin{thm}[\cite{AGr,MRud}]\label{MRudThm}
\begin{align*}
    \GW_{0,\Delta}^{\log}(\AAA,\Psi) = \GW_{0,\Delta}^{\trop}(\AAA,\Psi).
\end{align*}
\end{thm}
For the tropical disk counts determining $\alpha(\pp;p)$ in \eqref{TropFrobEqnQuot2}, one has $u_{(V_{\out},E_{\out})}=-p$. In particular, when $p=0$, these can be viewed as tropical curve counts.  We can therefore apply Theorem \ref{MRudThm} to replace the tropical disk counts in \eqref{TropFrobEqnQuot2} with the corresponding log invariants, yielding:
\begin{lem}\label{GWToric}
\begin{align*}
     \alpha(\pp;0)=\sum_{n\in K_2^+} \sum_{\ww\in \s{W}_{\pp}(n)} a_{\ww}\frac{ \GW_{0,\Delta_{\ww,\pp}}^{\log}(\AAA_{\ww,\pp,Q},\Psi_{\ww,\pp})}{|\Aut(\ww)|}z^{n}.
\end{align*}
\end{lem}

\section{Degeneration}\label{SectionDegen}

In this section, we will use a degeneration of our cluster varieties to relate the toric log Gromov-Witten invariants from Lemma \ref{GWToric} to log Gromov-Witten invariants of the cluster variety.  First though, we will need an important new technical result which says that curves satisfying ``somewhat generic'' conditions are torically transverse.

\subsection{Toric transversality lemma}\label{trans}

\begin{ntn}
We will write $\AAA_{\ww,\pp,Q}^{\circ}$ to indicate the tropical incidence conditions $\AAA_{\ww,\pp,Q}$ as in \S \ref{ThetaFromTrop}, except with each $\A_{ij}$ chosen to contain the origin in $\?{N}_{\bb{R}}$ (rather than being chosen to be a generic translate).  We write $\AAA_{\ww,\pp,0}^{\circ}$ to indicate that $Q=0$ as well. 
\end{ntn}

Consider the invariants $\GW_{0,\Delta_{\ww,\pp}}^{\log}(\AAA_{\ww,\pp,Q},\Psi_{\ww,\pp})$ of Lemma \ref{GWToric}, and consider a curve $\varphi^{\dagger}=[\varphi^{\dagger}:C^{\dagger}\rar \TV(\wt{\Sigma})^{\dagger}]\in \s{M}^{\log}_{0,\Delta_{\ww,\pp}}(\TV(\wt{\Sigma})^{\dagger},[\Delta_{\ww,\pp}])$. 
 We say that $\varphi^{\dagger}$ satisfies \textbf{somewhat generic} incidence and $\psi$-class conditions if it satisfies generically chosen representatives of the conditions $Z_{\AAA_{\ww,\pp,Q}}(\out)=:y^{\out}$ and $\Psi_{\ww,\pp}$, along with representatives of $Z_{\AAA_{\ww,\pp,Q}}((i,j))$ for each $(i,j)\in J_{\ww,\pp}$ which are not necessarily generic, but which at least intersect the interior of $D_i$.  The tropicalization of $\varphi^{\dagger}$ (in the sense of \cite{GSlog})\footnote{In this version of tropicalization from \cite[\S 1.4]{GSlog}, when a component of the basic log curve maps to the toric stratum corresponding to a cone $\sigma\in \wt{\Sigma}$, there is corresponding vertex of the tropical curve which can live anywhere in $\sigma$, with the precise location in $\sigma$ depending on a choice of pullback to the standard log point.  Nodes (respectively, marked points) of the log curve then correspond to compact edges (respectively, non-compact edges) of the tropicalization.  The somewhat-genericness of the incidence conditions ensures that their tropicalizations pass through the origin.\label{TropFoot}} is in $\f{T}_{0,\Delta_{\ww,\pp}}(\AAA_{\ww,\pp,0}^{\circ},\Psi_{\ww,\pp})$.  
 We wish to prove the following:\footnote{Note that Lemma \ref{DaggerCircLem} makes sense (and will hold) even in the $s=1$ case thanks to our convention of including the extra marked point $x_{\infty}$.}
 
\begin{lem}[Toric transversality lemma]\label{DaggerCircLem}
Suppose $[\varphi^{\dagger}:C^{\dagger}\rar \TV(\wt{\Sigma})^{\dagger}]\in \s{M}^{\log}_{0,\Delta_{\ww,\pp}}(\TV(\wt{\Sigma})^{\dagger},[\Delta_{\ww,\pp}])$ satisfies somewhat generic incidence and $\psi$-class conditions representing the classes $Z_{\AAA_{\ww,\pp,Q}}$ and $\Psi_{\ww,\pp}$ as in the invariants of Lemma \ref{GWToric}.  Assume $\Delta_{\ww,\pp}$ is a torically transverse degree for $\wt{\Sigma}$.  Then $\varphi(C)$ is torically transverse in $\TV(\wt{\Sigma})$.
\end{lem}

First, we need the following simple observation:
\begin{lem}\label{StarTransverse}
A basic stable log curve in $\TV(\wt{\Sigma})^{\dagger}$ is torically transverse if and only if all of its tropicalizations are supported on the $1$-skeleton of $\wt{\Sigma}$.
\end{lem}

Let $\star(0)\subset \?{N}_{\bb{R}}$ denote the union of all rational-slope lines through the origin.  Let $[\varphi^{\dagger}:C^{\dagger}\rar \TV(\wt{\Sigma})^{\dagger}]\in \s{M}^{\log}_{0,\Delta_{\ww,\pp}}(\TV(\wt{\Sigma})^{\dagger},[\Delta_{\ww,\pp}])$ be an arbitrary basic stable log map which satisfies our somewhat generic incidence and $\psi$-class conditions. By the assumption that $\Delta_{\ww,\pp}$ is a torically transverse degree for $\wt{\Sigma}$, the tropicalization $(\Gamma,h)$ of $\varphi^{\dagger}$ (for any pullback of the log structure to the standard log point) is supported on the $1$-skeleton of $\wt{\Sigma}$ if and only if it is supported on $\star(0)$.  Thus, to prove Lemma \ref{DaggerCircLem}, it suffices to prove that any $(\Gamma,h)$ obtained as a tropicalization of $\varphi^{\dagger}$ must be supported on $\star(0)$.

We will need some new definitions regarding tropical curves.  By a {\bf contractible} tropical curve we will mean a tropical curve as before, but now we allow compact positive-weight edges to be contracted by $h$.  Each flag $E\ni V$ is still assigned a designated primitive direction $u_{(V,E)}$ (i.e., as part of the data of the contractible tropical curve) such that the balancing condition still holds and such that $u_{(V,E)}=-u_{(V',E)}$ for $V,V'$ the two vertices of $E$.   
 One can define the type of a contractible tropical curve just as for the tropical curves of \S \ref{TropGW}, keeping in mind that contracted edges have directions.  Tropical curves as in \S \ref{TropGW} (i.e., without contracted compact positive-weight edges) will sometimes be referred to as \textbf{contracted} tropical curves.

Given a contractible tropical curve $(\Gamma',h')$, we can obtain a contracted tropical curve as follows: for each compact positive-weight edge of $\Gamma'$ contracted by $h'$, we simply contract the edge in the domain before applying $h'$ to get a new domain $\Gamma$.  Then $h$ is the map $\Gamma\rar N_{\bb{R}}$ such that the contraction $\Gamma'\rar \Gamma$ composed with $h$ is equal to $h'$.  We call this new tropical curve $(\Gamma,h)$ the {\bf contraction} of $(\Gamma',h')$, and we say that $(\Gamma',h')$ is an {\bf expansion} of $\Gamma$.  We say that a contractible tropical curve is in some $\f{T}_{0,\Delta}(\AAA,\Psi)$ if its contraction is.

The point of this terminology is that the moduli space of basic stable log maps is stratified by tropical types (cf. \cite{ACGSdecomp}), with the type of a stratum $\sigma$ being a contraction of the type associated to any stratum of $\partial \sigma$.

Consider $\Gamma\in \f{T}_{0,\Delta}(\AAA,\Psi)$, possibly contractible.  As in the computation of $\Mult(\Gamma)$ in \eqref{MultDfn}, we choose a flow on $\Gamma$ by specifying a vertex $V_0$ to serve as the sink. We recursively define affine linear spaces $\A_E$ associated to each edge $E\in \Gamma^{[1]}$ as follows.  For each $j\in J$, we take $\A_{E_j}:=\A_j$. For each vertex $V\neq V_0$, if $E_1,\ldots,E_k$ are the edges flowing into $V$ and $E_{V}$ is the unique edge flowing out of $V$, then we define
\begin{align}\label{AEdef}
    \A_E:=\bb{R}u_{E_V}+\bigcap_{i=1}^k \A_{E_i}.
\end{align}
Then for each $E\ni V_0$, we define $\A_{V_0,E}:=\A_E$.  Of course, we could take $V_0$ to be any vertex of $\Gamma^{[0]}$, and in this way we obtain linear spaces $\A_{V,E}$ for all flags of $\Gamma$.  One sees by induction that each vertex $V$ must be contained in $\A_{V,E}$ for each $E\ni V$.

Let us now specialize to the case of $\f{T}_{0,\Delta_{\ww,\pp}}(\AAA_{\ww,\pp,Q}^{\circ},\Psi_{\ww,\pp})$.  Recall from \S \ref{nontoric} that the seed data $\SSS$ included a form $B$ on $N$ which is skew-symmetrizable in the sense that there exists a skew-symmetric form $\omega$ on $N$ and positive integers $\{d_i\}_{i\in I_{\uf}}$ such that $B(e_i,e_j)=d_i\omega(e_i,e_j)$ for all $i,j\in I_{\uf}$.  In particular, this implies that $\ker(\omega|_{N_{\uf}})=\ker(\pi_2|_{N_{\uf}})$, so $\omega$ induces a non-degenerate skew-symmetric form $\?{\omega}$ on $\?{N}_{\uf}:=\pi_2(N_{\uf})$.  Furthermore, we see that $\pi_1(e_i)|_{\?{N}_{\uf}} = d_i\omega(e_i,\cdot)|_{\?{N}_{\uf}}=d_i\?{\omega}(\pi_2(e_i),\cdot)$ for each $i\in I_{\uf}$.  Hence, the conditions $\A_{ij}$ of $\AAA_{\ww,\pp,Q}^{\circ}$ satisfy $\A_{ij}\cap \?{N}_{\uf,\bb{R}} = u_{E_{ij}}^{\?{\omega}\perp}$, where for $u\in \?{N}_{\uf}$, $u^{\?{\omega}\perp}:=\{n\in \?{N}_{\uf,\bb{R}}:\?{\omega}(u,n)=0\}$.  This motivates the following:

\begin{lem}\label{AE}
Let $\Gamma\in \f{T}_{0,\Delta_{\ww,\pp}}(\AAA_{\ww,\pp,Q}^{\circ},\Psi_{\ww,\pp})$.  For $V\in \Gamma^{[0]}$ and $E\ni V$,  Let $\Gamma_{(V,E)}$ denote the closure in $\Gamma$ of the connected component of $\Gamma\setminus V$ which contains the interior of $E$.  Suppose all unbounded edges of $\Gamma_{(V,E)}$ are labelled by pairs $(i,j)\in J_{\ww,\pp}$.  Then $u_E\in \?{N}_{\uf}$ and $\A_{V,E}\cap \?{N}_{\uf,\bb{R}}\subset u_E^{\?{\omega}\perp}$.
\end{lem}
\begin{proof}
This follows from induction:  The claim holds for the edges $E_{ij}$ by the observations preceding the lemma.  If $E_1,\ldots,E_k$ are the edges flowing into a vertex $V'$ and $E'$ is the edge flowing out, then $u_{E_i}\in \?{N}_{\uf}$ for each $i=1,\ldots,k$ implies the same for $u_{E'}$ by the balancing condition.  Similarly, $\A_{E_i}\subset u_{E_i}^{\?{\omega}\perp}$ for each $i$ implies that $$\bigcap_{i=1}^k \A_{E_i} \subset \left(\sum_{i=1}^k \ww(E_i)u_{E_i}\right)^{\?{\omega}\perp}=u_{E'}^{\?{\omega}\perp},$$
hence
\begin{align*}
    \bb{R}u_{E'}+\bigcap_{i=1}^k \A_{E_i} \subset u_{E'}^{\?{\omega}\perp}, 
\end{align*}
as desired.
\end{proof}

The following lemma says that $Q$ being generic is enough to ensure that the tropical curves in $\f{T}_{0,\Delta_{\ww,\pp}}(\AAA_{\ww,\pp,Q}^{\circ},\Psi_{\ww,\pp})$ resemble $s$-tuples of broken lines meeting at a point.

\begin{lem}\label{Vout}
Fix a generic $Q\in N_{\bb{Q}}$.  Let $(\Gamma,h)\in \f{T}_{0,\Delta_{\ww,\pp}}(\AAA_{\ww,\pp,Q}^{\circ},\Psi_{\ww,\pp})$.  Let $V_{\out}$ be the vertex of $\Gamma$ contained in $E_{\out}$.  Then each component of $\Gamma \setminus V_{\out}$ other than (the interior of) $E_{\out}$ contains precisely one unbounded edge of the form $E_k$ for $k=1,\ldots,s,\infty\in J_{\ww,\pp}$.
\end{lem}
\begin{proof}
Suppose a component $\Gamma'$ does not contain any such edge.  Let $E'$ be the edge of $\Gamma'$ whose closure in $\Gamma$ contains $V_{\out}$.  Then by Lemma \ref{AE}, $\AAA_{V_{\out},E'}$ has codimension at least one, and so the generic point $Q$ (and thus $V_{\out}$)  cannot be contained in $\AAA_{V_{\out},E'}$.  This gives a contradiction.

On the other hand, $\Psi_{\ww,\pp}$ ensures that the valence of $V_{\out}$ is at least $(s+2)$, so no component of $\Gamma\setminus V_{\out}$ can contain more than one of the edges from $\{E_1,\ldots,E_s,E_{\infty},E_{\out}\}$.
\end{proof}

\begin{myproof}[Proof of Lemma \ref{DaggerCircLem}] 
A choice of generic $Q\in \?{N}_{\bb{Q}}$ determines (after a finite base change) a deformation $y_t^{\out}$ of the point $y^{\out}=:y_0^{\out}$ into the boundary of $\TV(\wt{\Sigma})$, cf. \cite[\S 3.2.1]{MRud}.  Denote the limit in the boundary by $y_1^{\out}$.  Since $y^{\out}$ was chosen generically, any curve $[\varphi^{\dagger}:C^{\dagger}\rar \TV(\wt{\Sigma})^{\dagger}]$ satisfying the somewhat generic conditions with $y^{\out}=y_0^{\out}$ will deform to a log curve $\varphi^{\dagger}_t$ satisfying the conditions with $y^{\out}$ replaced by $y_t^{\out}$.  The curves $\varphi_0^{\dagger}$ and $\varphi_1^{\dagger}$ admit tropicalizations $(\Gamma,h)\in \f{T}_{0,\Delta_{\ww,\pp}}(\AAA_{\ww,\pp,0}^{\circ},\Psi_{\ww,\pp})$ and $(\Gamma_1,h_1)\in \f{T}_{0,\Delta_{\ww,\pp}}(\AAA_{\ww,\pp,Q}^{\circ},\Psi_{\ww,\pp})$, respectively.  Furthermore, since $\varphi_1^{\dagger}$ is a degeneration of $\varphi_0^{\dagger}$, $(\Gamma_1,h_1)$ must be the same type as some expansion $(\Gamma',h')$ of $(\Gamma,h)$.

Suppose $\Gamma$ is not supported on $\star(0)$.  The incidence conditions at least force $h(V_{\out})=0$, where $V_{\out}$ denotes the vertex contained in $E_{\out}$.  Let $V'$ be a vertex of $\Gamma'$ of minimal distance from $V_{\out}$ which is not at $0$ and whose adjacent edges have directions not all parallel to the ray through $h(V')$.  For $E'\ni V'$ the edge on the component of $\Gamma'\setminus V'$ containing $V_{\out}$, we necessarily have $u_{E'}$ parallel to the ray through $h(V')$, and so $\A_{V',E'}$ contains the line $\ell_{V'}$ through $0$ and $h(V')$.  By Lemma \ref{Vout}, the component of $\Gamma'\setminus V_{\out}$ containing $V'$ includes exactly one of the edges $E_k$, $k=1,\ldots,s$.  It follows from this along with \eqref{AEdef}, Lemma \ref{AE}, and the genericness of $Q$ that $\A_{V',E'}$ has dimension at most $2$.
 
 If the dimension is $2$, then Lemma \ref{AE} applies to the other edges $E_i$, $i=1,\ldots,l$ containing $V'$, so they satisfy $\A_{V',E_i}\cap \?{N}_{\uf,\bb{R}}\subset u_{E_i}^{\?{\omega}\perp}$.  Since balancing forces these edges to point in multiple directions and $\?{\omega}$ is non-degenerate, these codimension $1$ spaces are non-equal, hence intersect to give a space of codimension at least $2$.  But since each $\A_{V',E_i}$ contains $V'$ and $0$, they must contain $\ell_{V'}$, so the intersection with $\A_{V',E'}$ is $1$-dimensional, hence non-transverse.  However (giving the edges of $\Gamma_1$ the same names as the corresponding edges of $\Gamma'$), translating $Q$ moves $h_1(E')$ independently from the other edges $E_i$ containing $V'$.  So $h_1(V')$ is contained in a translate of the space $\A_{V',E'}$ which has dimension at most $2$, but also in the space $\bigcap_{i=1}^{l} \A_{V',E_i}$, which has codimension at least $2$ and intersects $\A_{V',E'}$ non-transversely, thus giving a contradiction.

Similarly, if $\A_{V',E'}$ is $1$-dimensional, then still one of the $\A_{V',E_i}$ is contained in $u_{E_i}^{\omega\perp}$, which again must contain $\ell_{V'}$, hence have $1$-dimensional intersection with $\A_{V',E'}$.  This again is impossible for generic translates of $Q$ by the same argument.  The claim follows.
\end{myproof}

\subsection{Relating the invariants of $\TV(\wt{\Sigma})^{\dagger}$ and $(Y_{\SSS,\wt{\Sigma}},D_{\SSS,\wt{\Sigma}})$.}\label{DegenS}
Recall from Remark \ref{BlowupRmk} that the pair $(Y_{\SSS,\wt{\Sigma}},D_{\SSS,\wt{\Sigma}})$ can be constructed by, for each $i\in I_{\uf}$, blowing up the scheme-theoretic intersection
\begin{align*}
\?{H}_i:=D_{\pi_2(e_i)}\cap Z((a_i+z^{\pi_1(e_i)})^{|\pi_2(e_i)|}) \subset \TV_{\?{N}}(\wt{\Sigma})    
\end{align*}
for some $a_i\in \kk^*$, possibly followed by some toric blowdowns which (by \cite{AW}) do not affect log invariants. Lemma \ref{GWToric} allows us to express the theta functions associated to $(Y_{\SSS,\wt{\Sigma}},D_{\SSS,\wt{\Sigma}})$ in terms of certain log Gromov-Witten numbers of $\TV_{\?{N}}(\wt{\Sigma})$.  We wish to relate these to log Gromov-Witten numbers of $(Y_{\SSS,\wt{\Sigma}},D_{\SSS,\wt{\Sigma}})$.  More precisely, we wish to prove the following:

\begin{prop}\label{DegenProp}
Given a tuple $\pp=(p_1,\ldots,p_s)$ of vectors in $\?{N}$, consider the tropical degree $\Delta_{\pp}$ as in \eqref{Deltapp}.  
 Given $\beta\in \NE(Y_{\SSS,\wt{\Sigma}})$, let $\s{W}(\beta)$ denote the set of weight vectors $\ww$ such that $\sum_{j=1}^{l_i} w_{ij} = \beta.[E_i]$.  Let $Z(s+1)=[\pt]$ and $Z(i)=[Y_{\SSS,\wt{\Sigma}}]$ for all other $i$.  Let $\Psi(s+1)=s-1$ and $\Psi(i)=0$ for all other $i$.  Let $\eta$ denote the toric blowdown $Y^{\dagger}_{\SSS,\wt{\Sigma}}\rar Y^{\dagger}_{\SSS,\Sigma}$.  Then
\begin{align}\label{DegenFormula}
    \GW^{\log}_{0,Y_{\SSS,\wt{\Sigma}}^{\dagger},\Delta_{\pp},\beta}(Z,\Psi)z^{\eta_*\beta} = \sum_{n\in K_2^+} \sum_{\ww\in \s{W}_{\pp}(\varphi(p)+n)} a_{\ww}\frac{\GW^{\log}_{0,\Delta_{\ww,\pp}}(\AAA_{\ww,\pp,Q},\Psi_{\ww,\pp})}{|\Aut(\ww)|}z^{n}.
\end{align}
\end{prop}
Here, $z^{n}$ is viewed as an element of $\kk\llb \NE(Y_{\SSS,\wt{\Sigma}})\rrb$ by using $\kappa$ to identify $K_2$ with $N_1(Y_{\SSS,\Sigma})$.  Remark \ref{Dints} tells us already that $\kappa(n)$ is indeed equal to $\eta_*(\beta)$.

The strategy is to take a degeneration of $Y_{\SSS,\wt{\Sigma}}$, pictured in Figure \ref{CDfig}, as is done for two-dimensional cases in \cite[\S 5]{GPS}. 
To do this, let $\?{N}':=\?{N}\oplus \bb{Z}$, and let $\wt{\Sigma}^{\times}$ denote the fan in $N'_{\bb{R}}$ equal to the product of $\wt{\Sigma}$ with the fan for $\bb{A}^1$.   That is, for each cone $\sigma\in \wt{\Sigma}$, there are two cones in $\wt{\Sigma}^{\times}$ given by $\sigma\times \{0\}$ and $\wt{\sigma}:=\sigma \times \bb{R}_{\geq 0}$.  Then $\TV(\wt{\Sigma}^{\times})=\TV(\wt{\Sigma})\times \bb{A}^1$, and the projection $t:\?{N}'\rar \bb{Z}$ induces a map of fans giving a projection $t^{\times}:\TV(\wt{\Sigma}^{\times}) \rar \bb{A}^1$.

For any $n\in \?{N}$, let $\rho_n$ denote the ray generated by $n$ in $\?{N}$ (or the origin if $n=0$), and recall that for $n\neq 0$, $n'$ denotes the primitive vector in $\?{N}$ with direction $n$.  Now for each $i\in I_{\uf}$, we refine the cone $\wt{\rho}_{\pi_2(e_i)}\subset \wt{\Sigma}^{\times}$ by adding the ray $\bb{R}_{\geq 0}(\pi_2(e_i)',1)$.  We then further refine the cones of $\wt{\Sigma}^{\times}$ until we achieve a non-singular fan $\wt{\Sigma}'$ such that the projection $t':\TV(\wt{\Sigma}')\rar \bb{A}^1$ induced by $t$ is projective.  

Consider the cones in $\wt{\Sigma}'$ which are contained in $\wt{\rho}_{\pi_2(e_i)}$ for some $i\in I_{\uf}$, or in $\wt{\rho}_{p_i}$ for some $p_i$ from $\pp$.  Let $\wt{\Sigma}^{\circ}$ denote the set of all such cones except for the ones which are entirely contained in $(\?{N}_{\bb{R}},0)\subset \?{N}'_{\bb{R}}$.  Let $\?{Y}_0^{\circ}\subset \TV(\wt{\Sigma}')$ denote the union of all toric strata of $\TV(\wt{\Sigma}')$ which correspond to cones in $\wt{\Sigma}^{\circ}$.  Let $\?{D}_0^{\circ}\subset \?{Y}_0^{\circ}$ denote the union of all codimension-$2$ toric strata of $\TV(\wt{\Sigma}')$ which correspond to $2$-dimensional cones in $\wt{\Sigma}^{\circ}$.  By assuming that $\wt{\Sigma}$ was sufficiently refined, we can assume that $\?{D}_0^{\circ}$ is non-singular (in particular, the top-dimensional strata of $\?{D}_0^{\circ}$ are disjoint).

Finally, we blow up $\TV(\wt{\Sigma}')$ along the subvariety cut out by the loci $\?{H}_i':=Z(a_i+z^{(\pi_1(e_i),0)})\cap D_{(\pi_2(e_i),0)}$ for each $i\in I_{\uf}$. 
 Let $\wt{E}_i$, $i\in I_{\uf}$, denote the respective exceptional divisors.  We denote the resulting projective log smooth family by 
\[
\pi^{\dagger}:\wt{Y}_{\SSS,\wt{\Sigma}}^{\dagger} \rar (\bb{A}^1)^{\dagger},
\]
where the log structure of $\wt{Y}_{\SSS,\wt{\Sigma}}^{\dagger}$ is the divisorial log structures with respect to the proper transform of the toric boundary of $\TV(\wt{\Sigma}')$, and for $(\bb{A}^1)^{\dagger}$ we use the divisorial log structure with respect to $0\in \bb{A}^1$. 

Let $Y_0^{\circ}$ denote the preimage of $\?{Y}_0^{\circ}$ under the blowups of the loci $\?{H}_i$.  Similarly, let $D_0^{\circ}$ denote the proper transform of $\?{D}_0^{\circ}$ under these blowups.  Let $\wt{Y}_{\SSS,\wt{\Sigma}}^{\circ}$ denote the space $\wt{Y}_{\SSS,\wt{\Sigma}}$ equipped with the divisorial log structure associated to the divisor $Y_0^{\circ}$.  We then equip $Y_0^{\circ}$ with the log structure pulled back via the inclusion of $Y_0^{\circ}$ into $\wt{Y}_{\SSS,\wt{\Sigma}}^{\circ}$.  For each $i\in I_{\uf}$, let $\Bl_i$ denote the strata of $Y_0^{\circ}$ corresponding to cones in $\wt{\rho}_{\pi_2(e_i)}$.   Similarly, for each $p_i$, $i=1,\ldots,s$, let $\Bl_{p_i}$ denote the strata of $Y_0^{\circ}$ corresponding to cones in $\wt{\rho}_{p_i}$.

Let $Y_t^{\dagger}$ denote $\pi^{-1}(t)$ with the log structure induced by the inclusion into $\wt{Y}_{\SSS,\wt{\Sigma}}^{\dagger}$.  Note that for $t\neq 0$, $Y_t^{\dagger}$ is simply the cluster variety $Y_{\SSS,\wt{\Sigma}}^{\dagger}$ whose log GW invariants we are interested in, while $Y_{0}^{\dagger}$ includes 
\begin{align}\label{Y0}
Y_0^{\circ}=\TV(\wt{\Sigma})\cup \bigcup_{i\in I_{\uf}} \Bl_i
\end{align}
(with slightly different log structure) along with some additional strata.  Let $E_{i,t}:= \wt{E}_i \cap Y_{t}^{\dagger}$.  Let $\Bl'_i$ denote the component of $\Bl_i$ containing $E_{i,0}$, and let $\?{\Bl}'_i$ denote the image of $\Bl'_i$ under blowing down the $E_{j,0}$'s which it contains.  Let $\?{F}_i:=Z(a_i+z^{(\pi_1(e_i),0)})\cap \?{\Bl}'_i$, and let $F_i$ be the proper transform of $\?{F}_i$ under the blowups of the loci $\?{H}'_j$.  There is a fibration of $\?{\Bl}'_i$ with generic fibers (those not contained in the boundary of $\?{\Bl}'_i$) being $\bb{P}^1$, and with the two components of $D_0^{\circ}\cap \?{\Bl}_{E_{i,0}}$ being sections.  $\?{F}_i$ can be viewed as the union of the fibers which intersect $\?{H}'_i$.     Let $C_{F_i}\in \NE(\Bl'_i)_{\bb{Q}}$ be the class of the proper transform of one of these $\bb{P}^1$-fibers in $\?{F}_i$.  Let $D^{\circ}_{0,i}$ be the component of $D^{\circ}_0$ which intersects $F_i$.

Since $\pi^{\dagger}$ is log smooth, \cite[Thm. A.3]{MRud} says that the log Gromov-Witten invariants do not depend on $t$.  Hence, when proving Proposition \ref{DegenProp}, we can replace $(Y_{\SSS,\wt{\Sigma}},D_{\SSS,\wt{\Sigma}}) \cong Y_{t}^{\dagger}$ ($t\neq 0$) with $Y_{0}^{\dagger}$.

Consider the blowdown map $b:Y_0^{\dagger}\rar \TV(\wt{\Sigma})$.  Suppose $[\varphi^{\dagger}:C^{\dagger}\rar Y_0^{\dagger}]\in \s{M}^{\log}_{0,\Delta_{\pp}}(Y^{\dagger}_0,\beta)$ satisfies generic representatives of the incidence condition $\ev_{\out}^*([\pt])$ and the $\psi$-class conditions $\psi_{\out}^{s-1}$.  Then $b\circ \varphi:C\rar \TV(\wt{\Sigma})^{\dagger}$ can be equipped with a log structure making it into a curve in some $\s{M}^{\log}_{0,\Delta_{\ww,\pp}}(\TV(\wt{\Sigma})^{\dagger},[\Delta_{\ww,\pp}])$ which satisfies somewhat generic representatives of the conditions $Z_{\AAA_{\ww,\pp,Q}}$ and $\Psi_{\ww,\pp}$.  In particular, by Lemma \ref{DaggerCircLem}, $b\circ \varphi(C)$ must be torically transverse.  Hence $\varphi(C)$ must be torically transverse, and furthermore, any components of $\varphi(C)$ mapping to any $\Bl_i$ must be supported on fibers of the blowdown map $b$.  

We illustrate this setup in Figure \ref{CDfig}.

\begin{figure}[htb]
\noindent\begin{minipage}{0.45\textwidth}
\setlength{\parindent}{15pt}
    \caption{A sketch of the degeneration $\wt{Y}_{\SSS,\wt{\Sigma}}^{\dagger}$ of $Y_{\SSS,\wt{\Sigma}}^{\dagger}$.  The top part is a general fiber $Y_{\SSS,\wt{\Sigma}}$, while the bottom is the $0$-fiber $Y_0^{\dagger}$.  The bold part of $Y_0^{\dagger}$ forms $Y_0^{\circ}$, with the bolder beaded lines indicating $D_0^{\circ}$.  The thin curve in $\TV(\wt{\Sigma})$ represents a curve $b\circ \varphi(C)$. \label{CDfig}}
\end{minipage}
\hfill
\begin{minipage}{.5\textwidth}
\def\svgwidth{200pt}
    \input{Cluster_Degen.tex}
\end{minipage}
\end{figure}

\begin{lem}
Consider the locus of $\varphi^{\dagger}\in \s{M}^{\log}_{0,\Delta_{\pp}}(Y_0^{\dagger},\beta)$ satisfying generic representatives of the conditions $Z$ and $\Psi$.  Then the obstruction theory on this locus is unchanged if we view the maps $\varphi^{\dagger}$ as basic stable log maps to $Y_0^{\circ}$ instead of to $Y_0^{\dagger}$.
\end{lem}
\begin{proof}
We saw above that the image of such $\varphi^{\dagger}$ is necessarily torically transverse, so the only boundary divisors such curves can intersect are those in $D_0^{\circ}$.  Thus, forgetting the log structure along the other boundary divisors does not affect the obstruction theory.
\end{proof}

The upshot is that since $D_0^{\circ}$ is smooth, we can now use the log degeneration formula\footnote{Alternatively, we may use the main result of \cite{AMW} to say that the log Gromov-Witten invariants agree with the corresponding relative Gromov-Witten invariants, and then we may apply the relative degeneration formula of \cite{Li}.  Or as another alternative, one could use the recent log Gromov-Witten degeneration formula from expansions of \cite{RanExp} (allowing for non-smooth relative divisor), or a different log degeneration formula being developed in \cite{ACGS} which uses punctured invariants.}  of \cite[Thm. 1.4]{KLR}.   Let $[\varphi_{\TV}^{\dagger}:C^{\dagger}\rar \TV(\wt{\Sigma})^{\dagger}]\in \s{M}^{\log}_{0,\Delta_{\ww,\pp}}(\TV(\wt{\Sigma})^{\dagger},[\Delta_{\ww,\pp}])$ be a curve satisfying somewhat generic incidence and $\psi$-class conditions representing $\AAA_{\ww,\pp,Q}$ and $\Psi_{\ww,\pp}$.  The Gromov-Witten count of these curves is 
\begin{align}\label{TVGW}
    \GW^{\log}_{0,\Delta_{\ww,\pp}}(\AAA_{\ww,\pp,Q},\Psi_{\ww,\pp}).
\end{align} 
Furthermore, any curve $\wt{\varphi}^{\dagger}$ contributing to $\GW^{\log}_{0,Y_{\SSS,\wt{\Sigma}}^{\dagger},\Delta_{\pp},\beta}(Z,\Psi)$ is obtained by taking one of these curves $\varphi_{\TV}^{\dagger}$ and gluing chains of $\bb{P}^1$'s as follows:

For each $(i,j)$, we glue to the marked point $x_{ij}$ a chain of $b_i$ $\bb{P}^1$'s where $b_i\geq 1$ is the number of components of $\Bl_i$.  The first $(b_i-1)$ copies of $\bb{P}^1$ are just $w_{ij}|\pi_2(e_i)|$-fold covers of fibers of successive components of $\Bl_i$ with maximal tangency at $0$ and $\infty$ (where it intersects $D_0^{\circ}$).  The final $\bb{P}^1$, denoted $C_{ij}$, maps to $\Bl'_i$, satisfies $$(\wt{\varphi}|_{C_{ij}})_*[C_{ij}]=w_{ij}|\pi_2(e_i)|[C_{F_i}],$$ and has maximal tangency with $D_{0,i}^{\circ}$ at a point $p_{ij}$ ($p_{ij}$ is a node in $\wt{C}^{\dagger}$ but can be viewed as a marked point in $C_{ij}$ using the degeneration formula).   The nodes of this chain all have weight $w_{ij}|\pi_2(e_i)|$, and so the number of choices of log structures at these nodes, \textit{modulo automorphisms}, is 
\begin{align}\label{ijweight}
w_{ij}|\pi_2(e_i)|.
\end{align}
Since $Z_{\AAA_{\ww,\pp,Q}}$ imposes the condition that $\varphi_{\TV}$ maps $x_{ij}$ to $H_i$, we should (via a K\"unneth decomposition of the diagonal class) impose on $C_{ij}$ a condition that $p_{ij}$ maps to some curve $F_i^{\vee}\subset D^{\circ}_{0,i}$ which has intersection multiplicity $1$ with $F_i$ in $\Bl'_i$.  The resulting Gromov-Witten contribution of $C_{ij}$ is then reduced to the computation from \cite[Prop. 5.2]{GPS} (which was based on \cite[Thm. 5.1]{BP}), yielding  
\begin{align}\label{ijcover}
    \frac{(-1)^{w_{ij}-1}}{|\pi_2(e_i)|w_{ij}^2}.
\end{align}
Also, for each $p_i$, $i=1,\ldots,s$, we must glue a chain of copies of $\bb{P}^1$ in $\Bl_{p_i}$, each being a $|p_i|$-fold cover of a $\bb{P}^1$ in a fiber of $\Bl_{p_i}$ with maximal tangency at each intersection with $D_0^{\circ}$.  Such chains contribute a factor of $1$ to the Gromov-Witten count.

Finally, multiplying the toric Gromov-Witten count from \eqref{TVGW}, the node-weights $\prod_{ij} w_{ij}|\pi_2(e_i)|$ from \eqref{ijweight}, and the multiple-cover contributions $\prod_{ij} \frac{(-1)^{w_{ij}-1}}{|\pi_2(e_i)|w_{ij}^2}$ from \eqref{ijcover}, and then dividing by $|\Aut(\ww)|$ to correct for over-counting caused by labellings of the $x_{ij}$'s that are no longer remembered, we obtain by the degeneration formula that the contribution to \eqref{DegenFormula} of the curves coming from degree $\Delta_{\ww,\pp}$ is precisely the corresponding term from the right-hand side of \eqref{DegenFormula}, and then summing over all $\ww$ yields the desired result.

\qed

\begin{myproof}[Proof of Theorem \ref{MainThm}]
It is immediate from Proposition \ref{DegenProp} that $\Tr^s(\vartheta_{p_1},\ldots,\vartheta_{p_s})$, as defined in \eqref{Trs}, is indeed given as in \eqref{spointfunction}.  The theta functions were constructed as elements of a commutative associative algebra with $\vartheta_0=1$, so we already know that these properties are satisfied.  The fact that $\Tr^2$ and $\Tr^3$ uniquely determine the multiplication was Lemma \ref{nondegen}.  The identification of the base ring with $\kk\llb \NE(Y_{\SSS})\rrb$ is Lemma \ref{NECoeff}, and the finiteness statement for cases where the boundary supports an ample divisor was Corollary \ref{AmplePolynomial}.  The relation to the \cite{GHKK} theta functions on the Langlands dual cluster variety was Remark \ref{LangRmk}.
\end{myproof}

\section{The Gromov-Witten numbers are naive counts}\label{naiveSection}
Here we show that the log Gromov-Witten numbers $\GW^{\log}_{0,Y_{\SSS,\wt{\Sigma}}^{\dagger},\Delta_{\pp},\beta}(Z,\Psi)$ of Proposition \ref{DegenProp} are in fact naive counts of rational curves, not just virtual counts (assuming interior-curve freeness).  We denote by 
\begin{align*}
    \forget:\s{M}^{\log}_{0,\Delta}(Y_{\SSS,\wt{\Sigma}}^{\dagger},\beta)\rar \?{\s{M}}_{0,s+2}
\end{align*}
the forgetful/stabilization morphism taking $[\varphi^{\dagger}:C^{\dagger}\rar Y_{\SSS,\wt{\Sigma}}^{\dagger}]$ to the stabilization of the marked curve $C^{\dagger}$.  We allow any $s\geq 1$.

\begin{prop}\label{naive}
Suppose $[\varphi^{\dagger}:C^{\dagger}\rar Y_{\SSS,\wt{\Sigma}}^{\dagger}]\in  \s{M}^{\log}_{0,\Delta_{\pp}}(Y_{\SSS,\wt{\Sigma}}^{\dagger},\beta)$ satisfies generically chosen representatives for the incidence and $\psi$-class conditions $Z$ and $\Psi$ as in Proposition \ref{DegenProp}.  Then $\varphi(C)$ is torically transverse.  Let $[\pt]_Y$ denote the class of a point in $Y_{\SSS,\wt{\Sigma}}$ and let $[\pt]_{\s{M}}$ denote the class of a point in $\?{\s{M}}_{0,s+2}$.  Then the condition $\psi_{\out}^{s-1}$ can be replaced by $\forget^*[\pt]_{\s{M}}$, i.e.,
\begin{align}\label{psiforget}
\GW^{\log}_{0,Y_{\SSS,\wt{\Sigma}}^{\dagger},\Delta_{\pp},\beta}(Z,\Psi) =  \int_{[\s{M}^{\log}_{0,\Delta}(Y_{\SSS,\wt{\Sigma}}^{\dagger},\beta)]^{\vir}} \forget^*[\pt]_{\s{M}} \cup \ev_{\out}^*[\pt]_Y.
\end{align}
Furthermore, if $(Y_{\SSS,\wt{\Sigma}},D_{\SSS,\wt{\Sigma}})$ is interior-curve free, then \eqref{psiforget} is given by the naive count of irreducible torically transverse genus $0$ algebraic curves of tropical degree $\Delta_{\pp}$ and class $\beta$, with generically specified image under $\forget$, and with marked point $x_{\out}$ mapping to a generically specified point $y\in Y_{\SSS,\wt{\Sigma}}$.  These are the only basic stable log maps satisfying the generically specified representatives of the point and $\psi$-class conditions.
\end{prop}
\begin{proof}
For the toric transversality statement, we recall that by Lemma \ref{DaggerCircLem}, we had toric transversality in the central fiber of the degeneration of \S \ref{DegenS}, so the claim here follows from the fact that toric transversality is an open condition.

To prove \eqref{psiforget}, let $\?{\psi}_{\out}$ denote the corresponding $\psi$-class on $\?{\s{M}}_{0,s+2}$.  It is standard that $\?{\psi}_{\out}^{s-1}$ is the class of a point in $\?{\s{M}}_{0,s+2}$.  Furthermore, $\psi_{\out}-\forget^*\?{\psi}_{\out}$ is supported on the locus where the forgetful map destabilizes the curve-component $C_{\out}$ containing $x_{\out}$.  If $C_{\out}$ is destabilized when forgetting the map, then it intersects $D_{\SSS,\wt{\Sigma}}$ in at most one point because all such intersections must be marked points or nodes by \cite[Rmk 1.9]{GSlog}.  But Lemma \ref{ExpectedDim} below ensures that an irreducible genus $0$ curve hitting $D_{\SSS,\wt{\Sigma}}$ in at most one point will not hit the generically specified point $y\in Y_{\SSS,\wt{\Sigma}}$, contradicting the requirement that $x_{\out}$ maps to $y$.

For the claim about naive counts, we will first show that any $C$ satisfying the described generically specified conditions is irreducible.  Suppose $C$ were reducible, and let $C'$ be a component not containing $x_{\out}$.  We proceed by induction on $s$.

Using the condition $\forget^* [\pt]_{\s{M}}$, we can assume $C$ contains no contracted components.  Combined with the fact that $C$ is torically transverse and the interior-curve free assumption, it follows that no components of $C$ can map entirely into the boundary.  Let $C_0$ denote the closure in $C$ of the component of $C\setminus C'$ which contains $x_{\out}$, and let $C_1$ denote the closure in $C$ of $C\setminus C_0$.   By the interior-curve free assumption, $C'$ intersects the boundary, and as before, \cite[Rmk 1.9]{GSlog} tells us that any intersection of $C$ with $D_{\SSS,\wt{\Sigma}}$ is either a node or a marked point.  It follows that $C_1$ contains at least one of the marked points $x_i$, $i=1,\ldots,s$.

By the interior-curve free assumption and the fact that $C$ has no contracted components, $C_0$ must also intersect the boundary.  This already provides a contradiction when $s=1$, thus proving our base-case.  Furthermore, it now follows from the $\forget^* [\pt]_{\s{M}}$ condition that $C_1$ cannot contain more than one marked point.  Hence, $C_1$ must be irreducible, and now Lemma \ref{ExpectedDim} applies as before to say that the image of $C_1$ must lie in a fixed locus $E$ of codimension at least one which we can assume does not contain $y$. For convenience, let us assume that the unique $x_i$ in $C_1$ is $x_s$.

Next, note that applying $\forget$ destabilizes $C_1$, and upon stabilization, $x_s$ is identified with the point $x_s':=C_0\cap C_1$.  Treating $x_s'$ as an interior marked point on $C_0$, one obtains a basic stable log map $\varphi^{\dagger}_0:C^{\dagger}_0\rar Y_{\SSS,\wt{\Sigma}}^{\dagger}$.  Let $\varphi'_0:C'_0\rar Y_{\SSS,\wt{\Sigma}}^{\dagger}$ be the basic stable log map obtained by forgetting the marking at $x_s'$.  By the inductive assumption, $C_0$ must be irreducible.  We can therefore apply Lemma \ref{ExpectedDim} to say that generically specifying $y$ and the underlying marked curve $(C_0,x_1,\ldots,x_{s-1},x_{\out},x_{\infty})$ is sufficient to determine $\varphi_0'$ up to finitely many choices.  But then the condition that $x_s'$ maps to $E$ further determines the location of $x_s'$ up to finitely many choices, even though this point should independently be generically specified by the $\forget^*[\pt]_{\s{M}}$ condition.  This contradiction completes the proof of the claim that $C$ is irreducible.

We have thus shown that all curves satisfying the imposed conditions must lie in the open substack of irreducible curves $\s{M}^{\log,\irr}_{0,\Delta_{\pp}}(Y_{\SSS,\wt{\Sigma}}^{\dagger},\beta)\subset \s{M}^{\log}_{0,\Delta_{\pp}}(Y_{\SSS,\wt{\Sigma}}^{\dagger},\beta).$  The result now follows from another application of Lemma \ref{ExpectedDim}.
\end{proof}

The author learned of the following result and proof from Sean Keel and Tony Yue Yu, who are using a similar argument in \cite{KY}.  A version of the argument has previously appeared in \cite[Proof of Prop. 5.1]{Yu2}.
\begin{lem}\label{ExpectedDim}
Let $\Delta:J\rar \?{N}$ be a tropical degree and let $\beta\in \NE(Y_{\SSS,\wt{\Sigma}})$.  Let $\s{M}^{\log,\irr}_{0,\Delta}(Y_{\SSS,\wt{\Sigma}}^{\dagger},\beta)\subset \s{M}^{\log}_{0,\Delta}(Y_{\SSS,\wt{\Sigma}}^{\dagger},\beta)$ denote the open substack parametrizing basic stable log maps with irreducible domain curve $C$.  Label one of the markings as $x_{\out}$, with $\ev_{\out}$ being the corresponding evaluation map.  Then there exists a codimension $1$ subscheme $V\subset Y_{\SSS,\wt{\Sigma}}$ such that $\s{M}_V^{\circ}:=\ev_{\out}^{-1}(Y_{\SSS,\wt{\Sigma}}\setminus V) \cap \s{M}^{\log,\irr}_{0,\Delta}(Y_{\SSS,\wt{\Sigma}}^{\dagger},\beta)$ is smooth and has the expected dimension, i.e., $\dim(\s{M}_V^{\circ})=\dim(Y_{\SSS,\wt{\Sigma}})+|J|-3$.
\end{lem}
\begin{proof}
By \cite[Ch. III, Prop. 10.6]{Hart}, the derivative $d\ev_{\out}$ is surjective except on some codimension $1$ subscheme of $Y_{\SSS,\wt{\Sigma}}$, and we take this subscheme to be $V$.  To prove the claim then, we wish to show that deformations of basic stable log maps in $\s{M}_V^{\circ}$ are unobstructed.  Let $T^{\log}Y^{\dagger}_{\SSS,\wt{\Sigma}}=TY_{\SSS,\wt{\Sigma}}(-D_{\SSS,\wt{\Sigma}})$ denote the log tangent bundle of  $Y^{\dagger}_{\SSS,\wt{\Sigma}}$.  Note that $d\ev_{\out}$ is given at a point $[\varphi^{\dagger}:C^{\dagger}\rar Y_{\SSS,\wt{\Sigma}}^{\dagger}]\in \s{M}^{\log,\irr}_{0,\Delta}(Y_{\SSS,\wt{\Sigma}}^{\dagger},\beta)$ by the restriction
\begin{align*}
    d\ev_{\out}:H^0(C,\varphi^*T^{\log}Y^{\dagger}_{\SSS,\wt{\Sigma}})\rar (\varphi^*T^{\log}Y^{\dagger}_{\SSS,\wt{\Sigma}})_{x_{\out}}.
\end{align*}
Since $C\cong \bb{P}^1$, we can apply \cite[Ch. II, Def.-Prop. 3.8]{Kollar} to say that $\varphi^*T^{\log}Y^{\dagger}_{\SSS,\wt{\Sigma}}$ is semi positive, i.e., is isomorphic to a direct sum of line bundles of the form $\s{O}(a_i)$ for various $a_i\geq 0$.  It follows that $H^1(C,\varphi^* T^{\log}Y^{\dagger}_{\SSS,\wt{\Sigma}})=0$, and so $\s{M}_V^{\circ}$ is unobstructed, as desired.
\end{proof}

\begin{myproof}[Proof of Theorem \ref{MainNaive}]
Proposition \ref{naive} says that we can indeed interpret the Gromov-Witten counts of Theorem \ref{MainThm} as the naive counts described in Theorem \ref{MainNaive}.  The result follows.
\end{myproof}

\appendix

\section{Relation to quantum cohomology}\label{appendix}

Here we explain how to view the structure from Conjecture \ref{FrobConj} (the Frobenius structure conjecture) as part of a natural extension of (small) quantum cohomology to the log setting.  We begin by describing quantum cohomology in a somewhat new way.

Let $Y$ be a smooth projective variety over our algebraically closed field $\kk$ of characteristic $0$.  Let $\NE(Y)$ be the cone of effective curve classes in $Y$ up to numerical equivalence, and define
\begin{align*}
    \QH^*(Y):=H^*(Y,\bb{Q}\llb \NE(Y)\rrb),
\end{align*}
viewed for now as a $\bb{Q}\llb \NE(Y)\rrb$-module.  We define a $\bb{Q}\llb \NE(Y)\rrb$-multilinear $s$-point function $$\langle \cdot \rangle:\QH^*(Y)^{\otimes s} \rar \QH^*(Y)$$ as follows.  Given $\alpha_1,\ldots,\alpha_s\in H^*(Y,\bb{Q})$, and letting $[Y]$ denote the Poincar\'e dual to the fundamental class of $Y$, define $$\langle \alpha_1,\ldots,\alpha_s\rangle\in H^*(Y,\bb{Q}\llb\NE(Y)\rrb)$$ in terms of Gromov-Witten invariants via:
\begin{align}\label{spointQH}
    \langle \alpha_1,\ldots,\alpha_s\rangle := \sum_{\beta\in \NE(Y)} z^{\beta} \int_{[\s{M}_{0,s+2}(Y,\beta)]^{\vir}} \ev_{y_1}^*(\alpha_1) \wedge \cdots \wedge \ev_{y_s}^*(\alpha_s) \wedge \psi_{y_{s+1}}^{s-1}\ev_{y_{s+1}}^*[Y] \wedge \ev_{s+2}^*[Y].
\end{align}
This is then extended $\bb{Q}\llb \NE(Y)\rrb$-multilinearly to all of $\QH^*(Y)^{\otimes s}$.

\begin{thm}\label{QHthm}
There exists a unique associative product $*$ on $\QH^*(Y)$ which makes this $\bb{Q}\llb \NE(Y)\rrb$-module into a $\bb{Q}\llb \NE(Y)\rrb$-algebra and which satisfies
   \begin{align}\label{qprod}
       \langle \alpha_1,\ldots,\alpha_s\rangle = \langle \alpha_1 \cdots \alpha_s\rangle
   \end{align}
   for all $s$ and all $\alpha_1,\ldots,\alpha_s\in \QH^*(Y)$.  Furthermore, $\QH^*(Y)$ with this product $*$ is the small quantum cohomology ring of $Y$.
\end{thm}
The proof is based on the axioms of Gromov-Witten theory and the Topological Recursion Relation.  We suggest \cite[pg. 39 and Prop. 2.12]{Gr} as a reference for these properties. 
\begin{proof}
For $s\geq 2$, we can apply the Fundamental Class Axiom of Gromov-Witten theory to rewrite the Gromov-Witten invariants from \eqref{spointQH} as
\begin{align*}
    \int_{[\s{M}_{0,s+1}(Y,\beta)]^{\vir}} \ev_{y_1}^*(\alpha_1) \wedge \cdots \wedge \ev_{y_s}^*(\alpha_s) \wedge \psi_{y_{s+1}}^{s-2}\ev_{s+1}^*[Y].
\end{align*}
Now for $s=2$, the Point Mapping Axiom implies that $\langle \alpha_1,\alpha_2\rangle$ equals the usual Poincar\'e pairing of $\alpha_1$ and $\alpha_2$.  For $s=3$, the Dilation Axiom implies that \begin{align*}
    \langle \alpha_1,\alpha_2,\alpha_3\rangle = \sum_{\beta\in \NE(Y)} z^{\beta} \int_{[\s{M}_{0,3}(Y,\beta)]^{\vir}} \ev_{y_1}^*(\alpha_1) \wedge \ev_{y_2}^*(\alpha_2) \wedge  \ev_{y_3}^*(\alpha_3),
\end{align*}
i.e., the usual $3$-point function.  So the requirement that $\langle \alpha_1*\alpha_2,\alpha_3\rangle = \langle \alpha_1,\alpha_2,\alpha_3\rangle$ is the usual defining property of the quantum cohomology product $*$.

It remains to check that the usual quantum cohomology product satisfies \eqref{qprod} for all $s$.  The $s=1$ case is trivial, and the $s=2$ case follows from the Point Mapping Axiom as above.  The general case then follows by inductively applying the Topological Recursion Relation.
\end{proof}

It now seems natural to wonder whether there is an analog $H_{\log}^*(Y^{\dagger},\bb{Z})$ in the log setting which yields a log quantum cohomology ring $\QH_{\log}^*(Y^{\dagger})$ via the same recipe as above.  Indeed, the prime fundamental classes of our \S \ref{Intro} are precisely the degree $0$ log Chow classes of \cite{Bar,Herr}, and one expects that the invariants of Definition \ref{Nbeta} could be defined as in \eqref{spointQH}, i.e.,
\begin{align*}
        \langle p_1,\ldots,p_s\rangle := \sum_{\beta\in \NE(Y)} z^{\beta} \int_{[\s{M}^{\log}_{0,s+2}(\wt{Y}^{\dagger},\beta)]^{\vir}} \ev_{y_1}^*(\alpha_{p_1}) \wedge \cdots \wedge \ev_{y_s}^*(\alpha_{p_s}) \wedge \psi_{y_{s+1}}^{s-1}\ev_{y_{s+1}}^*[\pt] \wedge \ev_{s+2}^*[Y].
\end{align*}
Here, $\alpha_{p_i}$ denotes the prime fundamental class associated to $p_i\in B(\bb{Z})$.  Also, note that restricting to degree $0$ here forces us to replace $\ev_{y_{s+1}}^*[Y]$ by $\ev_{y_{s+1}}^*[\pt]$ (the trace is only nonzero on top-degree elements, so we are using multiplication by the class of a point to pull the trace back to degree $0$ elements).  So this is the sense in which Conjecture \ref{FrobConj} is an extension of degree $0$ quantum cohomology to the log setting.

\begin{rmk}
In higher-degree, progress is obstructed by the fact that it is not even clear how to define $H_{\log}^*(Y^{\dagger},\bb{Z})$.  One idea is to use an extension of \cite{Bar,Herr}'s log Chow groups which includes ``punctured classes'' (conditions one can impose on punctured invariants), and the author hopes this will yield an appropriate definition of log quantum cohomology.

Alternatively, Ganatra and Pomerleano \cite{GP,GP2} suggest taking 
\begin{align*}
H_{\log}^*(Y^{\dagger},\bb{Z}):=\bigoplus_{p\in B(\bb{Z})} H^*(D_p^{\KN},\bb{Z}),    
\end{align*}
where $D_{p}^{\KN}$ denotes the Kato-Nakayama space \cite{KN} associated to $D_p^{\dagger}$ (i.e., $D_p$ with the log structure pulled back from inclusion into $Y^{\dagger}$).  Gross-Pomerleano-Siebert \cite{GPomS} seek to prove that working with punctured invariants and this choice of $H_{\log}^*(Y^{\dagger},\bb{Z})$ yields the symplectic cohomology ring $\SH_{\log}^*(Y^{\dagger})$ of $(Y,D)$.  Closed string mirror symmetry predicts that $\SH^*_{\log}(Y^{\dagger})$ is isomorphic to the ring of polyvector fields on the mirror, so in particular, $\SH^0_{\log}(Y^{\dagger})$ should be isomorphic to the coordinate ring of the mirror (cf. \cite[\S 1]{Pas}).    In this sense, the Frobenius structure conjecture should be the degree $0$ part of the closed string mirror symmetry conjecture, as noted in the original statement \cite[arXiv v1, Conj. 0.8]{GHK1}, and also in \cite[pg. 5]{GSInt}.  The observation that this construction is analogous to that of quantum cohomology was also previously noted by Gross and Siebert \cite[Rmk 2.3]{GSInt}.
\end{rmk}

\bibliographystyle{alpha}  
\bibliography{main}        
\index{Bibliography@\emph{Bibliography}}%

\end{document}